\documentclass[12pt]{article}

\usepackage{amsmath,amsthm}
\usepackage{amssymb}
\usepackage{epsfig}

\usepackage{easybmat}
\usepackage{graphicx,hyperref}

\newtheorem{theorem}{Theorem}[section]
\newtheorem{lemma}[theorem]{Lemma}

\theoremstyle{definition}
\newtheorem{define}[theorem]{Definition}
\newtheorem{remark}[theorem]{Remark}
\newtheorem{example}[theorem]{Example}

\newcommand{\dsty}{\displaystyle}

\newcommand{\pf}{\noindent {\bf Proof: }}
\newcommand{\enpf}{\hfill $\Box$ \vspace{.2in} }
\newcommand{\vs}{\vspace{0.15in}}
\newcommand{\nin}{\noindent}

\newcommand{\ord}[4]{({#1} \hspace{0.07in}  {#2} \hspace{0.07in}  {#3}  \hspace{0.07in}  {#4}) }

\begin{document}

\title{Existence and Stability of Four-Vortex Collinear \\ Relative 
Equilibria with Three Equal Vorticities}

\author{Brian Menezes\thanks{bjmene16@g.holycross.edu} \qquad Gareth E. Roberts\thanks{groberts@holycross.edu}  \\
  \\
Dept. of Mathematics and Computer Science \\
College of the Holy Cross}

\maketitle

\begin{abstract}
We study collinear relative equilibria of the planar four-vortex problem where three of the four vortex strengths are identical.
The $S_3$ invariance obtained from the equality of vorticities is used to reduce the defining equations and organize the solutions into two distinct
groups based on the ordering of the vortices along the line.  The number and type of solutions are given, along with a discussion of the bifurcations
that occur.  The linear stability of all solutions is investigated rigorously and stable solutions are found to exist for cases where the vorticities have mixed signs.  
We employ a combination of analysis and computational algebraic geometry to prove our results.
\end{abstract}


\vspace{.1in}


\section{Introduction}

The planar $n$-vortex problem is a Hamiltonian system describing the motion of $n$~point vortices in the plane acting
under a logarithmic potential function.    It is a well-known model for approximating vorticity evolution in fluid dynamics~\cite{aref-int, newton}. 
One of the most fruitful approaches to the problem is to study {\em stationary} configurations, solutions where the initial
configuration of vortices is maintained throughout the motion.  As explained by O'Neil~\cite{oneil}, there are four possibilities:
equilibria, relative equilibria (uniform rotations), rigidly translating configurations, and collapse configurations.  
Much attention has been given to relative equilibria since numerical simulations of certain physical processes (e.g., the eyewall
of hurricanes~\cite{davis, KossSchub}) often produce rigidly rotating configurations of vortices.   Analyzing the stability of relative equilibria
improves our understanding of the local behavior of the flow;  it also has some practical significance given the persistence of these
solutions in numerical models of hurricane eyewalls.   Other physical examples are provided in~\cite{aref-newton}.

There are many examples of stable relative equilibria in the planar $n$-vortex problem.  Perhaps the most well known is the equilateral 
triangle solution, where three vortices of arbitrary circulations are placed at the vertices of an equilateral triangle.  If the sum of the circulations
does not vanish, then the triangle rotates rigidly about the center of vorticity.  This periodic solution is linearly
(and nonlinearly) stable provided that the {\em total vortex angular momentum}
$
L \; = \;  \sum_{i < j}  \Gamma_i \Gamma_j 
$
is positive~\cite{synge}, where $\Gamma_i \in \mathbb{R} - \{0\}$ represents the circulation or vorticity of the $i$th vortex.  
Other stable examples include the regular $n$-gon for $4 \leq n \leq 7$ (equal-strength circulations required)~\cite{thomson, have, aref-equil, cs, schm};
the $1 + n$-gon for $n \geq 3$ (a regular $n$-gon with an additional vortex at the center)~\cite{cs};
the isosceles trapezoid~\cite{g:stability}; a family of rhombus configurations~\cite{g:stability};
and configurations with one ``dominant'' vortex and $n$ small vortices encircling the larger one~\cite{bhw, bhl}.
In~\cite{aref-stab3}, Aref provides a comprehensive study of three-vortex collinear relative equilibria, finding linearly
stable solutions for certain cases when the vortex strengths have mixed signs.
The rhombus configuration studied in~\cite{g:stability} and some particular solutions of the $(1+3)$-vortex problem discussed in~\cite{bhl} 
provide some additional examples of stable solutions with circulations of opposite signs.

Relative equilibria can be interpreted as critical points of the Hamiltonian $H$ restricted to a level surface of the angular impulse~$I$.  This gives
a promising topological viewpoint to approach the problem~\cite{palmore}.   If all vortices have the same sign, then a relative equilibrium is linearly stable
if and only if it is a nondegenerate minimum of $H$ restricted to $I = \mbox{constant}$~\cite{g:stability}.  Moreover, because $I$ is a conserved quantity,
a technique of Dirichlet's applies to show that any linearly stable relative equilibrium with same-signed circulations is also nonlinearly stable.

In this paper we apply methods from computational algebraic geometry to investigate 
the existence and stability of collinear relative equilibria in the four-vortex problem.  To make the problem
more tractable, we restrict to the case where three of the vortices are assumed to have the same circulation.  Specifically, if $\Gamma_i$ is
the circulation of the $i$th vortex, then we assume that $\Gamma_1 = \Gamma_2 = \Gamma_3 = 1$ and
$\Gamma_4 = m$, where $m \in \mathbb{R} - \{0\}$ is treated as a parameter.   Solutions to this problem come in groups of six due to the
invariance that arises from permuting the three equal-strength vortices.   We use this invariance to simplify the problem considerably, obtaining
a complete classification of the number and type of solutions in terms of~$m$.   We also provide a straight-forward algorithm to rigorously find
all solutions for a fixed $m$-value.

When counting the number of solutions, we follow the usual convention (inherited from the companion problem in celestial mechanics) 
of identifying solutions that are equivalent under rotation, scaling, or translation.  In other words, we count equivalence classes of relative
equilibria.  In general, there are $n!/2$ ways to arrange $n$ vortices on a common line, where the factor of $1/2$ occurs because configurations
equivalent under a $180^\circ$ rotation are identified.   The 12 possible orderings in our setting are organized into two groups.  Group~I contains
the 6 arrangements where the unequal vortex (vortex~4) is positioned exterior to the other three;  Group~II consists of the 6 orderings
where vortex~4 is located between two equal-strength vortices.  We show that for any $m > -1/2$, there are exactly 12 collinear relative equilibria, 
one for each possible ordering of the vortices.  As $m$ decreases through $-1/2$, the solutions in Group~II disappear;  there are 
precisely 6 solutions for each $m \in (-1, -1/2]$, one for each ordering in Group~I.  There are no collinear relative equilibria for $m \leq -1$.

We also consider the linear stability of the collinear relative equilibria in the planar setting.  Due to the integrals and symmetry that naturally arise
for any relative equilibrium, there are always four {\em trivial} eigenvalues $0, 0, \pm i$ (after a suitable scaling).  
For the case $n = 4$, there are four nontrivial eigenvalues remaining that determine stability.  We explain how the nontrivial eigenvalues
can be computed from the trace~$T$ and determinant~$D$ of a particular $2 \times 2$ matrix and provide useful formulas 
for $T$ and~$D$ as well as conditions that guarantee linear stability.   By applying these formulas and conditions to our
specific problem, we are able to rigorously analyze the linear stability of all solutions in Groups I and~II.

We show that the Group~II solutions are always unstable, with two real pairs of nontrivial eigenvalues $\pm \lambda_1, \pm \lambda_2$.
The Group~I solutions go through two bifurcations, at $m = m_c \approx -0.0175$ and $m = m^{\ast} \approx -0.8564$.  These important
parameter values are roots of a particular sixth-degree polynomial in~$m$ with integer coefficients.  For $m > m_c$, the Group~I solutions are unstable,
with two real pairs of eigenvalues.  At $m = m_c$, these pairs merge and then bifurcate into a complex quartuplet $\pm \alpha \pm i \beta$
for $m \in (m^{\ast}, m_c)$.  The Group~I solutions are linearly stable for $m \in (-1, m^{\ast})$ and spectrally stable at $m = m^{\ast}$.
The linear stability of the Group~I solutions is somewhat surprising since four of the six solutions limit on
a configuration with a pair of binary collisions as $m \rightarrow -1^{+}$.

This problem has recently been explored in~\cite{perez}, where the intent was to classify all relative equilibria, not just the collinear configurations.  
Unfortunately, there are some errors in this paper.  For example, Theorem~4 claims the existence of two families of rhombus configurations.
However, this violates the main theorem in~\cite{albouy}, which states that a convex relative equilibrium is symmetric with respect to one diagonal
if and only if the circulations of the vortices on the other diagonal are equal.   To obtain a rhombus, there must be two pairs of equal-strength vortices, one pair for
each diagonal  (see Section~7.4 in~\cite{HRS} for the complete solution).  If three vortices have equal circulations, then the only possible rhombus relative equilibrium 
is a square.   In this article we treat the collinear case in much greater depth than in~\cite{perez} and focus on the linear stability of solutions (the stability question is not
considered in~\cite{perez}).


Much of our work relies on the theory and computation of Gr\"{o}bner bases and 
would not be feasible without the assistance of symbolic computing software.
Computations were performed using Maple~\cite{maple} and 
many results were checked numerically with Matlab~\cite{matlab}.  
The award-winning text by Cox, Little, and O'Shea~\cite{CLO} is an excellent reference 
for the theory and techniques used in this paper involving modern and computational algebraic geometry.

The paper is organized as follows.  In the next section we introduce relative equilibria and provide the set up for our particular
family of collinear configurations.   We then explain how the solutions come in groups of six and use the invariance inherent
in the problem to locate, count, and classify solutions in terms of the parameter~$m$.   
In Section~3 we provide the relevant theory and techniques for studying the linear stability of relative equilibria in the planar $n$-vortex problem.
Applying these ideas in our specific setting, we obtain reductions that reduce the stability problem to the calculation
of two quantities, $T$ and~$D$.  This leads to the discovery of linearly stable solutions and the bifurcation values that signify
a change in eigenvalue structure.

\section{Collinear Relative Equilibria with Three Equal Vorticities}
\label{sec:rel-eq}

We begin with some essential background.  The planar $n$-vortex problem was first described as
a Hamiltonian system by Kirchhoff~\cite{kirchhoff}.
Let $z_i \in \mathbb{R}^2$ denote the position of the $i$th vortex and let
$r_{ij} = \|z_i - z_j\|$ represent the distance between the $i$th and $j$th vortices.
The mutual distances $r_{ij}$ are useful variables.  
The motion of the $i$th vortex is determined by
\begin{equation}
\Gamma_i  \dot z_i \;  = \;  J \frac{\partial H}{\partial z_i}  \; = \;   J  \sum_{j \neq i}^n \frac{\Gamma_i \Gamma_j}{r_{ij}^2}(z_j - z_i), \quad 1 \leq  i \leq n ,
\label{eq:vorts}
\end{equation}
where 
$
J=
\begin{bmatrix}
0 & 1\\
-1 & 0
\end{bmatrix}
$
is the standard $2 \times 2$ symplectic matrix and
$$
H = -\sum_{i<j} \Gamma_i \Gamma_j \ln (r_{ij})
$$
is the Hamiltonian function for the system.
The {\em total circulation} of the system is $\Gamma = \sum_i  \Gamma_i $, and as long as $\Gamma \neq 0$,
the {\em center of vorticity}  $c =  \frac{1} {\Gamma} \sum_i \Gamma_i z_i$ is well-defined.

\subsection{General facts about relative equilibria}

\begin{define}
A {\em relative equilibrium} is a periodic solution of~(\ref{eq:vorts}) where each vortex rotates about~$c$
with the same angular velocity $\omega \neq 0$.   Specifically, we have
$$
z_i(t) = c + e^{- \omega J t} (z_i(0) - c) , \quad \mbox{ for each } i \in \{1, \ldots , n\}.
$$ 
\end{define}
It is straight-forward to check that the mutual distances $r_{ij}$ in a relative equilibrium are unchanged throughout the motion, so that
the initial configuration of vortices is preserved.  
Upon substitution into system~(\ref{eq:vorts}), we see that the initial positions of a relative equilibrium must satisfy the following system
of algebraic equations
\begin{equation}
- \omega \, \Gamma_i (z_i  - c) \; = \;  \frac{\partial H}{\partial z_i}  \; = \;  \sum_{j \neq i}^n \frac{\Gamma_i \Gamma_j}{r_{ij}^2} (z_j - z_i), \quad  i \in \{1, \ldots, n\}.
\label{eq:rel-equ}
\end{equation}

Suppose that $z_0 = (z_1(0), \ldots, z_n(0))$ represents the initial positions of a relative equilibrium.   Although it is a periodic solution,
it is customary to treat a relative equilibrium as a point~$z_0 \in \mathbb{R}^{2n}$ (e.g., a fixed point in rotating coordinates).
We will adopt this approach here.
From equation~(\ref{eq:rel-equ}), we see that any translation, scaling, or rotation of~$z_0$ leads to another
relative equilibrium (perhaps with a different value of $c$ or~$\omega$).   Thus, relative equilibria are never isolated and it makes
sense to consider them as members of an equivalence class~$[[ z_0 ]]$,
where $w_0 \sim z_0$ provided that $w_0$ is obtained from $z_0$ by translation, scaling, or rotation.    
The stability type of $z_0$ is the same for all members of~$[[ z_0 ]]$.  Reflections of $z_0$ are also relative equilibria
(e.g., multiplying the first coordinate of $c$ and each $z_i$ by~$-1$), but
these will not be regarded as identical when counting solutions.

The quantity 
$$
I \;  = \;  \sum_{i=1}^n \, \Gamma_i \|z_i - c\|^2  \; = \;   \frac{1}{\Gamma}  \sum_{i < j}  \Gamma_i \Gamma_j r_{ij}^2
$$
can be regarded as a measure of the relative size of the system.  
It is known as the {\em angular impulse} with respect to the center of vorticity, the analog of the {\em moment of inertia} in
the $n$-body problem.  The angular impulse is an integral of motion for the planar $n$-vortex problem~\cite{newton}.
One important property of~$I$ is that relative equilibria (regarded as points in~$\mathbb{R}^{2n}$) are critical
points of the Hamiltonian restricted to a level surface of $I$.  This can be seen by rewriting 
system~(\ref{eq:rel-equ}) as
\begin{equation}
\nabla H(z) + \frac{\omega}{2} \nabla I(z) \; = \; 0 ,
\label{eq:cc}
\end{equation}
where $\nabla$ is the usual gradient operator.   Here we treat the constant $\omega/2$ as a Lagrange multiplier.
This gives a very useful topological approach to the study of relative equilibria.  The main result of~\cite{g:stability} is that,
for positive vorticities, a relative equilibrium is linearly stable if and only if it is a nondegenerate minimum of~$H$ restricted to $I = \mbox{constant}$.
Using equation~(\ref{eq:cc}), it is straight-forward to derive the formula $\omega = L/I$.

\subsection{Defining equations}

\begin{figure}[hbt]
\begin{center}
\includegraphics[width=400bp]{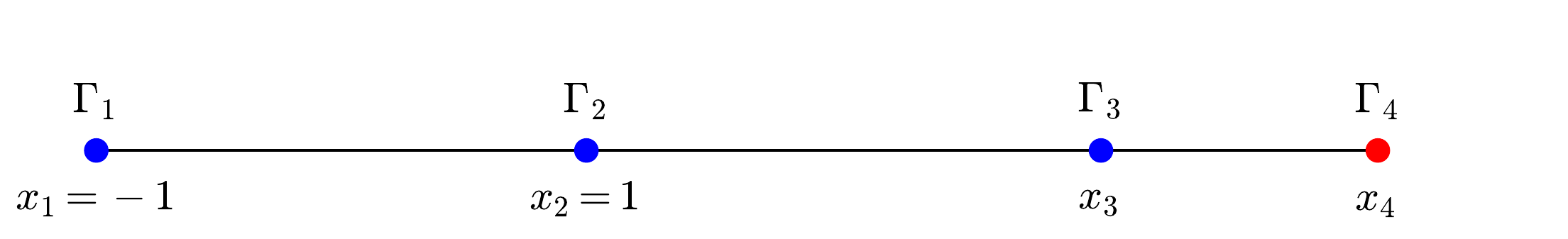}
\hspace{0.1in}
\includegraphics[width=400bp]{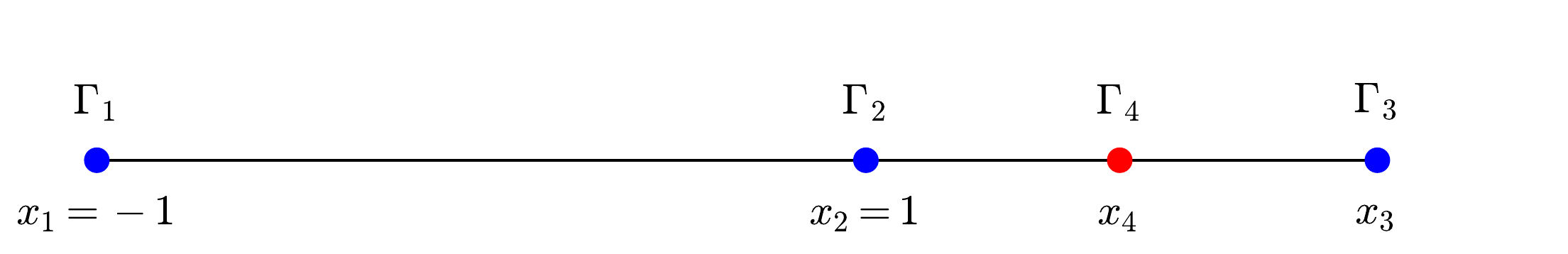}
\end{center}
\vspace{-0.2in}
\caption{Two collinear relative equilibria, each with circulations $\Gamma_1 = \Gamma_2 = \Gamma_3 = 1, \Gamma_4 = -0.25$.
The top configuration has ordering $\ord{1}{2}{3}{4}$, with $x_3 \approx 3.104$ and $x_4 \approx 4.228$, while the ordering of the bottom solution
is $\ord{1}{2}{4}{3}$, with $x_3 \approx 2.328$ and $x_4 \approx 1.659$.}
\label{fig:setup}
\end{figure}

We now focus on four-vortex relative equilibria whose configurations are collinear, that is, all vortices
lie on a common line.  To make the problem tractable, we assume that three of the four vortex strengths are identical.
Set $\Gamma_{1} = \Gamma_{2} = \Gamma_{3}=1$, and $\Gamma_{4}=m$, where $m \in \mathbb{R}$ is a parameter.
Without loss of generality, we take the positions of the relative equilibrium to be on the $x$-axis, $z_i = (x_i, 0)$, and translate
and scale the configuration so that $x_1 = -1$ and $x_2 = 1$.  This produces a simpler system to solve than other approaches
(e.g., setting $c = 0$ and $\omega = 1$).  It also helps elucidate the inherent symmetries in the problem. 
In our set up, the center of vorticity $c$ and angular vorticity $\omega$ will vary, 
but the coordinates of the first two vortices will remain fixed (see Figure~\ref{fig:setup}).

According to equation~(\ref{eq:rel-equ}), a relative equilibrium in this form must satisfy the following system of equations:
\begin{eqnarray}
\omega(-1 - c) + \frac{1}{2} + \frac{1}{x_3 + 1} +  \frac{m}{x_4 + 1} & = & 0,  \nonumber \\[0.07in]  \nonumber
\omega(1 - c) - \frac{1}{2} + \frac{1}{x_3 - 1} +  \frac{m}{x_4 - 1} & = & 0,  \\[0.07in]   \label{eq:orgSys}
\omega(x_3 - c) - \frac{1}{x_3 + 1} - \frac{1}{x_3 - 1} +  \frac{m}{x_4 - x_3} & = & 0,  \\[0.07in]  \nonumber
\omega(x_4 - c) - \frac{1}{x_4 + 1} - \frac{1}{x_4 - 1} -  \frac{1}{x_4 - x_3} & = & 0. 
\end{eqnarray}

Let the numerators of the left-hand side of each equation above be denoted by $f_1, f_2, f_3,$ and $f_4$ respectively, and append the two polynomials
$$
f_5 \; = \; u(x_4 - 1) - 1  \quad \mbox{ and } \quad  f_6 \; = \; v(x_4 + 1) - 1,
$$
in order to eliminate solutions with collisions (i.e., $x_4 = \pm 1$).  Let $F$ be the polynomial ideal generated by $f_1, \ldots, f_6$
in $\mathbb{Q}[\omega, c, u, v, x_4, x_3, m]$.  Computing a Gr\"{o}bner basis of $F$, denoted by $GB$, with respect to the lexicographic order
$\omega > c > u > v > x_4 > x_3 > m$, yields a basis with 26 elements. The first of these is
a 12th-degree polynomial in $x_3$ with coefficients in $m$, given by
$$
\begin{array}{c}
P(x_3, m) = (m+1)(2m+1)(m+2)^2 x_3^{12} - (32m^5 + 224m^4 + 635m^3 + 873m^2 + 576m + 144)x_3^{10}  \\[0.05in]
+ (640m^5+4066m^4+10126m^3+12546m^2+7776m+1944)x_3^8\\[0.05in]
- (3776m^5 + 23984m^4 + 60278m^3  + 75042m^2 + 46656m +11664)x_3^6  \\[0.05in]
+ (5760m^5+40806m^4+115191m^3+158841m^2+104976m+26244)x_3^4 \\[0.05in]
-27m^2(96m^3+464m^2+717m+351)x_3^2 + 54m^4.
\end{array}
$$
Although $P$ appears intimidating to analyze, we will use the equality of the vorticities and invariant group theory to factor it into
four cubic polynomials in $x_3$.

Before performing this reduction, we repeatedly apply the Extension Theorem to insure that a zero of $P$ can be extended
to a full solution of system~(\ref{eq:orgSys}).  The third term in $GB$ is 
$$
Q(x_4, x_3, m) \; = \;  192m^3(4m + 5)(x_3^2 + 3)x_4 \; + \;  q_1(x_3, m),
$$
where $q_1$ is a polynomial in the variables $x_3$ and $m$.  If $m \neq 0$ and $m \neq -5/4$, then the Extension Theorem applies
to extend a zero of $P$, call it $(x_3^\ast, m^\ast)$, to a solution $(x_4^\ast, x_3^\ast, m^\ast)$.  Moreover, since $Q$ is linear
in $x_4$, there is a unique such extension.   The 20th term in $GB$ is 
$$
(m+1)(c(m+3) - x_3 - mx_4),
$$
which implies that
$$
c \; = \;  \frac{x_3 + mx_4}{m+3}
$$
as long as $m \neq -1, -3$.  As expected, this agrees with the formula for the center of vorticity in our set up.  Similar arguments
work to extend any zero of $P$ uniquely to a solution of the full system~(\ref{eq:orgSys}).  Note that we have not ruled out the case
that $x_3^\ast = x_4^\ast$, a collision between the third and fourth vortices.  We have proven the following lemma.

\begin{lemma}
Fix an $m \in \mathbb{R}$ with $m \neq -3, -5/4, -1, 0$.  Then any solution $x_3^\ast$ to $P = 0$ can be extended uniquely to
a full solution of system~(\ref{eq:orgSys}).  
\label{lemma:extend}
\end{lemma}

\subsection{Classifying solutions}

Since vortices 1, 2, and 3 have the same vorticity, we can interchange their positions (a relabeling of the vortices) to create a new relative equilibrium.
However, because $x_1 = -1$ and $x_2 = 1$ are always assumed, it is necessary to apply a scaling and translation in order
to convert a relabeled solution into our specific coordinate system.  This creates a map between solutions
of system~(\ref{eq:orgSys}).  To make these ideas precise, we will keep track of how the vortices are
arranged under different permutations.  

\begin{define}
If vortices $i, j, k,$ and $l$ are positioned so that $x_i < x_j < x_k < x_l$, then the corresponding {\em ordering} is denoted
$\ord{i}{j}{k}{l}$.
\end{define}

To illustrate the inherent invariance in system~(\ref{eq:orgSys}), suppose that we have a relative equilibrium 
with coordinates 
$$
x \; = \;  (x_1, x_2, x_3, x_4)  \; = \;  (-1, 1, a, b), \quad  \mbox{ with } 1 < a < b.
$$  
This corresponds to the ordering $\ord{1}{2}{3}{4}$.  Relabeling the
vortices in the order $\ord{3}{1}{2}{4}$ gives another relative equilibrium, but with coordinates
\begin{equation}
x' \; = \;  (x_1', x_2', x_3', x_4')  \; = \;  (1, a, -1, b),
\label{eq:newRE}
\end{equation}
which does not match our setup.  The linear map
$\phi(x_i) = \frac{2 x_i - a - 1}{a - 1}$ satisfies $\phi(1) = -1$ and $\phi(a) = 1$.  Consequently, applying $\phi$ 
to each entry in~(\ref{eq:newRE}) will convert $x'$ into the correct form.  Because $\phi$ is a scaling and translation, the resulting 
coordinate vector,
$$
\phi(x') \; = \; \left( -1, 1, \frac{3+a}{1 - a}, \frac{-2b+1+a}{1 - a} \right),
$$ 
is also a relative equilibrium, and its first two coordinates match our setup.
This gives a new solution to system~(\ref{eq:orgSys}), one with ordering $\ord{3}{1}{2}{4}$.

The above argument demonstrates an important invariance for the ideal~$F$.
It shows that for a fixed value of $m$, if $(x_3, x_4) = (a,b)$ is a partial solution in the variety of $F$, then so is $(x_3, x_4) = (\frac{3+a}{1 - a}, \frac{-2b+1+a}{1 - a})$.  
Put another way, if $\overline{F} = F \cap \mathbb{Q}[x_3, x_4, m]$ is an elimination ideal, then $\overline{F}$ is invariant under the map
\begin{equation}
S(x_3, x_4) \; = \;  \left( \frac{3 + x_3}{1 - x_3}, \frac{-2x_4 + 1 + x_3}{1 - x_3} \right),
\label{eq:symmS}
\end{equation}
after clearing denominators.  (Here we can assume that $x_3 \neq 1$ because $x_3 = 1$ is a collision between
vortices 2 and~3.)  Another symmetry, which is easy to discern, arises by reflecting all four positions about the origin:
\begin{equation}
R(x_3, x_4) \; = \;  (-x_3, -x_4).
\label{eq:symmR}
\end{equation}
However, this operation reverses the ordering of the vortices (e.g., ordering $\ord{1}{2}{3}{4}$ maps to $\ord{4}{3}{2}{1}$), and thus requires
that vortices 1 and~2 be interchanged in order to insure that $x_1 = -1$ and $x_2 = 1$ is maintained (e.g., ordering $\ord{4}{3}{2}{1}$ becomes $\ord{4}{3}{1}{2}$).  
We can also consider the composition of $R$ and $S$ to generate additional invariants for $\overline{F}$.  As expected,
this yields a total of six invariants (including the identity) for $\overline{F}$, as $R$ and $S$ generate a group of order six that is isomorphic to $S_3$,
the symmetric group on three symbols.   

\begin{theorem}
Let $G$ be the group generated by the maps $R$ and $S$ under composition, where $R$ and $S$ are given by (\ref{eq:symmR}) and~(\ref{eq:symmS}), respectively.
Then $G$ is isomorphic to $S_3$ and the elimination ideal $\overline{F} = F \cap \mathbb{Q}[x_3, x_4, m]$ arising from system~(\ref{eq:orgSys})
is invariant under~$G$.   Consequently, solutions to system~(\ref{eq:orgSys}) come in groups of six.
\label{ThmInv}
\end{theorem}

\pf
As explained above, $\overline{F}$ is invariant under both $R$ and $S$.  This was confirmed using Maple by checking that, 
for each polynomial $p$ in a Gr\"{o}bner basis of $\overline{F}$,
$p(R(x_3, x_4))$ and the numerator of $p(S(x_3, x_4))$ are also in $\overline{F}$.  
It follows that $\overline{F}$ is also invariant under any composition of these maps.

Let $e$ represent the identity function $e(x_3, x_4) = (x_3, x_4)$.
We compute that $R^2 = R \circ R = e, S^3 = S \circ S \circ S = e$, and $(R \circ S)^2 = e$.  This is sufficient to show that $G$ is
isomorphic to $S_3$.

Since $\overline{F}$ is invariant under $G$ and the order of~$G$ is six, it follows that one solution in the variety of $\overline{F}$ leads to five others.  
The only possible exception occurs when $G$ has fixed points, that is, points in $\mathbb{R}^2$ that are mapped to the same place under different group transformations.  
A straight-forward calculation reveals that the only possible fixed points are $(0,0), (3,1), (-3,-1), (-1,b)$, and $(1,b)$, where $b \in \mathbb{R}$ is arbitrary.
However, each of these corresponds to a collision between two vortices, and is thus excluded.  Therefore, the six solutions generated by~$G$ are distinct.
\enpf

Suppose that we have a relative equilibrium solution
with $(x_3, x_4) = (a, b)$, where $1 < a < b$.  This solution has ordering $\ord{1}{2}{3}{4}$.  Applying the transformations from~$G$ generates
five additional solutions, each with a different ordering of the vortices.
These solutions and their corresponding orderings are shown in the first two columns of Table~\ref{table:orders} and will be denoted as {\em Group~I}.
Likewise, for a solution with $1 < b < a$, which corresponds to the ordering $\ord{1}{2}{4}{3}$, there are five other solutions generated by $G$ whose orderings are displayed
in the third column of Table~\ref{table:orders}.   These solutions will be referred to as {\em Group~II}.
The 12 orderings from the union of the two groups are the only allowable orderings because we have assumed that $x_1 < x_2$, thereby eliminating
half of the 24 permutations in~$S_4$.   Note that each of the orderings in Group~I have vortex 4 positioned exterior to the three equal-strength vortices, while
for Group~II, the fourth vortex always lies between two of the equal-strength vortices.

\renewcommand{\arraystretch}{2.3}
\begin{table}[h]
\begin{center}
\begin{tabular}{||c|c|c|c||}
\hline
\makebox[1.1in] {{\bf Coordinates} $(x_3, x_4)$} &  \makebox[1.1in]{{\bf Group I}} & \makebox[1.1in]{{\bf Group II}} &  \makebox[1.5in]{{\bf Group Element}} \\
\hline\hline
$(a,b)$ &  $\ord{1}{2}{3}{4}$ &  $\ord{1}{2}{4}{3}$ &  $e$   \\
\hline
$(-a, -b)$  & $\ord{4}{3}{1}{2}$  & $\ord{3}{4}{1}{2}$ &  $R$     \\
\hline
$\dsty{ \left(\frac{a+3}{a-1}, \frac{-2b+a+1}{a-1} \right) }$  & $\ord{4}{1}{2}{3}$  & $\ord{1}{4}{2}{3}$ &  $R \circ S$     \\
\hline
$\dsty{ \left(\frac{3 - a}{a+1}, \frac{2b-a+1}{a+1} \right) }$  & $\ord{1}{3}{2}{4}$  & $\ord{1}{3}{4}{2}$ &  $S \circ R$     \\
\hline
$\dsty{ \left(\frac{3+a}{1 - a}, \frac{-2b+a+1}{1 - a} \right) }$  & $\ord{3}{1}{2}{4}$  & $\ord{3}{1}{4}{2}$ &  $S$     \\
\hline
$\dsty{ \left(\frac{a-3}{a+1}, \frac{-2b+a-1}{a+1} \right) }$  & $\ord{4}{1}{3}{2}$  & $\ord{1}{4}{3}{2}$ &  $S^2 = S \circ S$     \\
\hline\hline
\end{tabular}
\caption{{\bf Invariants:} Since $\Gamma_1 = \Gamma_2 = \Gamma_3$, relative equilibria $x = (-1, 1, x_3, x_4)$ come in groups of six, each with a different ordering of vortices.  
The positions of the third and fourth vortices (first column) are expressed in terms of one particular solution $(x_3 = a, x_4 = b)$.  
The orderings in Group~I arise from a solution where $1 < a < b$, while those in Group~II satisfy $1 < b < a$.  The corresponding group
element from $G \simeq S_3$ is shown in the last column.  }
\label{table:orders}
\end{center}
\end{table}
\renewcommand{\arraystretch}{1.0}

\begin{remark}
Recall that $S_3$ is isomorphic to $D_3$, the dihedral group of degree three.  Since $G \simeq S_3 \simeq D_3$, the transformations
$S$ and $S^2$ from Table~\ref{table:orders} correspond to rotations, while the remaining non-identity elements represent reflections.
\end{remark}

\subsection{Using invariant group theory to find solutions}

Based on the discussion in the previous section, we can apply invariant group theory to rigorously study the solutions to system~(\ref{eq:orgSys})
in terms of the parameter~$m$.  Let $r_1 = x_3, r_2 = (x_3-3)/(x_3+1),$ and $r_3 = (3+x_3)/(1-x_3)$ denote the three values of $x_3$ corresponding to the group elements
$e, S^2,$ and $S$, respectively (the rotations).  The cubic polynomial with these three roots should be a factor of $P(x_3, m)$, the first polynomial in
the lex Gr\"{o}bner basis $GB$ arising from system~(\ref{eq:orgSys}).

We introduce the coordinates $\sigma, \tau,$ and $\rho$, defined by the elementary symmetric functions on the roots $r_1, r_2,$ and $r_3$:
\begin{eqnarray}
\sigma & = &  r_1 + r_2 + r_3,  \label{eq:sigma} \\
\tau & = &  r_1 r_2  + r_1 r_3 + r_2 r_3, \label{eq:tau} \\
\rho & = &  r_1 r_2 r_3 .    \label{eq:rho}
\end{eqnarray}
Consider the ideal in $\mathbb{Q}[\sigma, \tau, \rho, x_3, m]$ generated by equations (\ref{eq:sigma}), (\ref{eq:tau}), and~(\ref{eq:rho}) (after clearing denominators) 
and the twelfth-degree polynomial $P(x_3, m)$.  Computing a lex
Gr\"{o}bner basis (denoted $GB^\ast$) with respect to the ordering $\tau < \sigma < x_3 < \rho < m$ yields a basis with four polynomials, the first two of which are
\begin{eqnarray*}
P_1 & = &  (m+1)(2m+1)(m+2)^2 \rho^4  - m^2(32m^3 +152m^2 + 239m + 117)\rho^2 +54m^4, \\
P_2 & = &  x_3^3 + \rho x_3^2  - 9x_3 - \rho.
\end{eqnarray*}
The polynomial $P_2$ is just the expanded version of equation~(\ref{eq:rho}). 
The fact that $P_1$ is even in $\rho$ is expected from the reflection symmetry~$R$.  
Note that $(\rho, x_3) \mapsto (-\rho, -x_3)$ is a symmetry for $P_2 = 0$.

This Gr\"{o}bner basis calculation effectively factors~$P$
into the product of four cubics in the form of~$P_2$.  Indeed, using Maple, we confirm that
$$
P \; = \;  (m+1)(2m+1)(m+2)^2  \prod_{i = 1}^4  (x_3^3 + \rho_i x_3^2  - 9x_3 - \rho_i),
$$
where $\rho_i, i \in \{1, 2, 3, 4\}$, are the four roots of $P_1$ (each a function of~$m$).

We now analyze the roots of the quartic $P_1$ as a function of~$m$.

\begin{lemma}
The roots $\rho_i(m)$ of $P_1$ satisfy the following properties:
\begin{itemize}
\item For $m > -1/2$ and $m \neq 0$, $P_1$ has four real roots.
\item  For $-1 < m \leq -1/2$, $P_1$ has exactly two real roots.
\item  At $m = 0$, $P_1$ has precisely one real root at zero of multiplicity four.
\item  At $m = -1/2$, $P_1$ reduces to a quadratic function with two real roots at $\pm \sqrt{3/7}$.
\item  For $m \leq -1$, $P_1$ has no real roots.
\end{itemize}
\label{Lemma:QuarticRho}
\end{lemma}

\pf
Introduce the variable $\xi = \rho^2$.  The results follow by treating $P_1$ as a quadratic function of~$\xi$ with coefficients in $\mathbb{Q}[m]$.  
The discriminant of $P_1(\xi)$ is given by
$$
\Delta_1 \; = \;  m^4 (4m+5)^2 (64m^4 + 448m^3 + 1153m^2 + 1278m + 513),
$$
which is clearly positive for $m > 0$.  Using Sturm's Theorem~\cite{sturm}, it is straight-forward to check that $\Delta_1 > 0$
for $-1 < m < 0$ as well.  Therefore, the roots of $P_1(\xi)$ are real for $m > -1$.  

If $m > -1/2$ and $m \neq 0$, the leading coefficient and constant term of $P_1(\xi)$ are positive, while the middle term has a negative coefficient.
Since the roots are real, Descartes' Rule of Signs shows that $P_1(\xi)$ has two positive roots, which implies that $P_1(\rho)$ has four real roots
of the form $\pm \rho_1(m), \pm \rho_2(m)$.

The leading coefficient of $P_1$ becomes negative for $-1 < m < -1/2$, while the middle coefficient flips sign at $m \approx -0.942.$
Thus, the sign pattern for the coefficients of $P_1(\xi)$ is either $- \, - \, +$ or $- \, + \, +$.  In either case there is just one sign change, so
Descartes' Rule implies that $P_1(\xi)$ has only one positive root.  Thus, for $-1 < m < -1/2$, 
$P_1(\rho)$ has precisely two real roots of the form $\pm \rho_1(m)$.

For $m < -1$, all three coefficients are positive so we have $P_1(\xi) > 0$ for $\xi \geq 0$.  Consequently, $P_1(\rho)$ has no real roots.  
The remaining facts listed for the specific $m$-values $ 0, -1/2,$ and $-1$ are easily confirmed.
\enpf

\begin{lemma}
The roots of $P_2(x_3) = x_3^3 + \rho x_3^2  - 9x_3 - \rho$ are real and distinct for any $\rho \in \mathbb{R}$.  If $r_1 = a$ is a root, then the other
two roots are given by
$$
r_2 \; = \;  \frac{a - 3}{a + 1} \quad \mbox{and} \quad r_3 \; = \;  \frac{3 + a}{1 - a} \; .
$$
Let $r_1 = a$ denote the largest root.
If $\rho > 0$, then the roots satisfy $1 < r_1 < 3, -1 < r_2 < 0,$ and $r_3 < -3$. 
If $\rho < 0$, then the roots satisfy $r_1 > 3, 0 < r_2 < 1,$ and $-3 < r_3 < -1$. 
\label{Lemma:Cubic}
\end{lemma}

\pf
The discriminant of $P_2$ with respect to $x_3$ is $4(27 + \rho^2)^2$, which is always positive.  Consequently, the roots of $P_2$ are always real.
If $a$ is a root of $P_2$, then we have $\rho = (a^3 - 9a)/(1 - a^2) = r_1 r_2 r_3$, as expected.  Then, it is straight-forward to check that $P_2$ factors as
$(x - r_1)(x - r_2)(x - r_3)$.  

Next we note that $P_2(1) = -8$ and $P_2(3) = 8 \rho$.  By the Intermediate Value Theorem, we have a root $a$ satisfying $1 < a < 3$ if $\rho > 0$,
or $a > 3$ if $\rho < 0$.  In the first case, we see that $-1 < r_2 < 0$ and $r_3 < -3$ by straight-forward algebra.  This also serves to show that $r_1 = a$ is
the largest root.  For the case $\rho < 0$ and $a > 3$, we have $0 < r_2 < 1$ and $-3 < r_3 < -1$, as desired.  
\enpf

Lemma~\ref{Lemma:Cubic} is important because it provides specific information on the location of the third vortex without having to work with
the complicated expressions that arise from Cardano's cubic formula.  Note that if $\rho$ is a complex number, then the roots of $P_2$ must also be complex.

\vs

\nin  {\bf Algorithm for computing solutions:}  Applying Lemma~\ref{Lemma:Cubic} and the reductions outlined above, we have the following algorithm for computing
the positions of all relative equilibria for a fixed value of~$m$.  In theory, the calculations are exact because they only require solving, in order, a quadratic, cubic, and linear equation.

\begin{itemize}
\item[{\bf 1.}]  Compute the real roots $\rho_i$ of the even quartic $P_1$.

\item[{\bf 2.}]  For each real value of $\rho_i$, substitute into the cubic $P_2$ and find the largest root to obtain~$x_3$.

\item[{\bf 3.}]  Substitute $x_3$ and $m$ into $Q$ and solve $Q = 0$ for $x_4$.  (Recall that $Q$ is linear in $x_4$.)

\item[{\bf 4.}]  Two additional solutions for $(x_3, x_4)$ are obtained by using the formulas in the bottom two rows of Table~\ref{table:orders}.
\end{itemize}

\begin{remark}

\begin{enumerate}
\item  Each choice of $\rho$ leads to three distinct solutions with different orderings.
By symmetry, using both $\rho$ and $-\rho$ yields six solutions that correspond to six orderings in a particular group (either Group~I or Group~II).
Thus, two positive roots of the quadratic $P_1(\xi)$ will generate 12 solutions,
while one positive root leads to six solutions.  If $P_1(\xi)$ has no positive roots, then there are no solutions.

\item  The remaining two polynomials in the Gr\"{o}bner basis $GB^\ast$ reveal some peculiar properties of solutions.  The third polynomial in $GB^\ast$ is simply
$\sigma + \rho$, which implies that the sum of the roots of $P_2$ is the negative of the product of the roots.  The remaining entry in $GB^\ast$ is just
$\tau + 9$, which reveals that the symmetric product of the three roots is always equal to $-9$.  These facts can also be verified by examining the
coefficients of $P_2$ and are apparently an artifact of our special choice of coordinates $x_1 = -1, x_2 = 1$.  
\end{enumerate}

\label{remark:GBstar}
\end{remark}

Next we demonstrate our algorithm for finding all relative equilibria solutions in two important cases.

\begin{example}
{\bf The case $m=1$}.

If all four vortices have the same strength $\Gamma_i = 1$, then the four roots of $P_1$ are $\rho_i = \pm (\sqrt{3} \pm \sqrt{2})$.  
Taking $\rho = \sqrt{3} - \sqrt{2}$,  the three roots of $P_2$ are
\begin{equation}
-\sqrt{3} - \sqrt{2} \; \approx \; -3.146, \quad 1 + \sqrt{2} - \sqrt{6} \; \approx \; -0.035, \quad \mbox{and} \quad -1 + \sqrt{2} + \sqrt{6} \; \approx \; 2.864 \, .
\label{eq:rootsP2}
\end{equation}
Notice that the sum and product of these roots equals $-\rho$ and $\rho$, respectively, in accordance with part 2 of Remark~\ref{remark:GBstar}.
Since $r_{12} = 2$, we expect the solution with ordering $\ord{1}{2}{3}{4}$ to have $r_{34} = 2$ by symmetry.
Choosing $x_3 = -1 + \sqrt{2} + \sqrt{6}$ and $x_4 = 1 + \sqrt{2} + \sqrt{6} \approx 4.864$ gives the desired solution.
After computing the corresponding values of $c$ and $\omega$, 
this solution was confirmed by substituting it into the Gr\"{o}bner basis $GB$ as well as into system~(\ref{eq:orgSys}).
The center of vorticity is $c = (\sqrt{2} + \sqrt{6})/2 \approx 1.932$ and the angular velocity is $\omega = 3/(6 + 2\sqrt{3}) \approx 0.317$.  
The other two roots of $P_2$ shown in equation~(\ref{eq:rootsP2}) yield solutions with orderings 
$\ord{3}{1}{2}{4}$ and $\ord{4}{1}{3}{2}$, with $x_4$ coordinates $\sqrt{3} + \sqrt{2}$ and $-1 + \sqrt{2} - \sqrt{6}$, respectively.
These solutions concur with those obtained by using the symmetry transformations indicated on the bottom two rows of Table~\ref{table:orders}.

If we choose $\rho = -\sqrt{3} - \sqrt{2}$ instead, then we obtain $x_3 =  1 + \sqrt{2} + \sqrt{6} \approx 4.864$ as the largest root of $P_2$.
Then $x_4 = -1 + \sqrt{2} + \sqrt{6} \approx 2.864$ gives the coordinate of the fourth vortex.  This solution corresponds to ordering $\ord{1}{2}{4}{3}$.  
The other two roots of $P_2$ lead to solutions with orderings $\ord{3}{1}{4}{2}$ and $\ord{1}{4}{3}{2}$.  

The remaining two values of $\rho$ lead to six other solutions corresponding to the orderings in rows 2, 3, and 4 of Table~\ref{table:orders}.
All 12 solutions are symmetric, with the distance between the first pair of vortices equal to the distance between the second pair
(i.e., for the ordering $\ord{i}{j}{k}{l}$, we have $r_{ij} = r_{kl}$).  The ratio of this distance over the distance between the inner pair of vortices
is always $(\sqrt{3}+\sqrt{2}-1)/2 \approx 1.073$.  All 12 solutions are geometrically equivalent.  These results agree with those given
in Section~5.1 of~\cite{HRS}.
\enpf
\label{ex:equal}
\end{example}

\begin{example}
{\bf The case $m=0$}.

The solutions for the case $\Gamma_4 = 0$ are relative equilibria of the restricted four-vortex problem.  The three equal-strength vortices are akin to the large masses (called {\em primaries}) in the
celestial mechanics setting.   It is known that the primaries must form a relative equilibrium on their own~\cite{xia}.
Thus, based on symmetry, we expect the first three vortices to be equally spaced.

Setting $m=0$ in $P_1$ gives $\rho_i = 0  \; \forall i$.  The cubic $P_2$ reduces to $x_3^3 - 9x_3$ and thus $P(x_3, m)$ factors as 
$$
P(x_3, 0)  \;  =  \;  4 x_3^4 (x_3 - 3)^4 (x_3 + 3)^4,
$$
with roots $x_3 = -3, 0, 3$ each repeated four times.   As expected, each of the three possible values for $x_3$ yield an equally-spaced
configuration for the equal-strength vortices.

The repeated roots and the fact that Lemma~\ref{lemma:extend} does not apply when $m=0$ suggest a bifurcation.
Surprisingly, this does {\em not} happen:  there are still 12 different
solutions, one for each possible ordering in Table~\ref{table:orders}.  
While the first nine elements of the Gr\"{o}bner basis $GB$, including $Q$, vanish entirely at $m=0$ and $x_3 = 3$,
the tenth element yields a quartic polynomial in $x_4$ with four distinct real roots.  The same feature occurs if $x_3 = -3$
or $x_3 = 0$.   We obtain 12 solutions given by
$$
(x_3, x_4) \; = \;  \left(  a_i, \, \frac{a_i}{3} \pm \frac{(a_i^2 + 9) \sqrt{54 \pm 6 \sqrt{57}} }{54} \right),
$$
where $a_i = -3, 0,$ or $3$, and all four sign combinations occur.  These solutions were checked to insure that each satisfied
system~(\ref{eq:orgSys}).

\enpf
\label{ex:zero}
\end{example}

We now have enough information to count and classify all solutions in terms of the parameter~$m$.

\begin{theorem}
The number and type of four-vortex collinear relative equilibria with circulations
$\Gamma_1 = \Gamma_2 = \Gamma_3 = 1$ and $\Gamma_4 = m$ are given as follows:

\begin{itemize}
\item[{\bf (i)}]  If $m > -1/2$, then there are 12 relative equilibria, one for each possible ordering of the vortices (both Groups I and~II are realized);
\item[{\bf (ii)}]  If $-1 < m \leq -1/2$, there are 6 relative equilibria, one for each ordering in Group~I;
\item[{\bf (iii)}]  If $m \leq -1$, there are no solutions.
\end{itemize}

\label{thm:main}
\end{theorem}

\pf
First, we compute the discriminant of $P(x_3, m)$ as a polynomial in~$x_3$, and find
that for $m > -1$, the discriminant vanishes only if $m = -1/2$ or $ m = 0$.  Consequently, the roots of $P$ are distinct for $m > -1$ except in these
two special cases.  

{\bf (i)}  Suppose that $m$ is fixed with $m > -1/2$ and $m \neq 0$.  Lemmas \ref{Lemma:QuarticRho} and~\ref{Lemma:Cubic} combine to show that $P$ has 12 distinct roots, 
and Lemma~\ref{lemma:extend} implies that each of these $x_3$-values can be extended to a full solution of system~(\ref{eq:orgSys}).  Thus there are 12 relative equilibria.
We now show that each possible ordering in Table~\ref{table:orders} is realized.  

Recall from Example~\ref{ex:equal} that for the case $m=1$, choosing $\rho = \sqrt{3} - \sqrt{2}$ (the smaller positive root of $P_1$) leads to a solution with $1 < x_3 < x_4$ (ordering $\ord{1}{2}{3}{4}$).
This solution can be continued analytically as $m$ varies away from 1 by following the solution corresponding to the
smaller positive root of $P_1$ and the largest root of $P_2$.   By continuity, the only way for the ordering $\ord{1}{2}{3}{4}$ to disappear is for there to be a collision of vortices,
either with $x_3 = x_2 = 1$ or $x_3 = x_4$, at some particular $m$-value.  The first of these possibilities is ruled out by Lemma~\ref{Lemma:Cubic}.  
The second possibility is eliminated by taking the third and tenth polynomials in $GB$ and making the substitution $x_4 = x_3$.  
Computing a Gr\"{o}bner basis for these two polynomials, along with $P_1$ and $P_2$, produces the polynomial~$1$.  
Consequently, there are no solutions in the variety of $F$ with $x_3 = x_4$.  We note that neither $x_3$ nor
$x_4$ can become infinite for a particular $m$-value when $m > -1/2$.  This follows from 
Lemmas \ref{Lemma:QuarticRho} and~\ref{Lemma:Cubic} and from the fact that $Q$ is linear in $x_4$.

\begin{figure}[htb]
\begin{center}
\includegraphics[width=220bp]{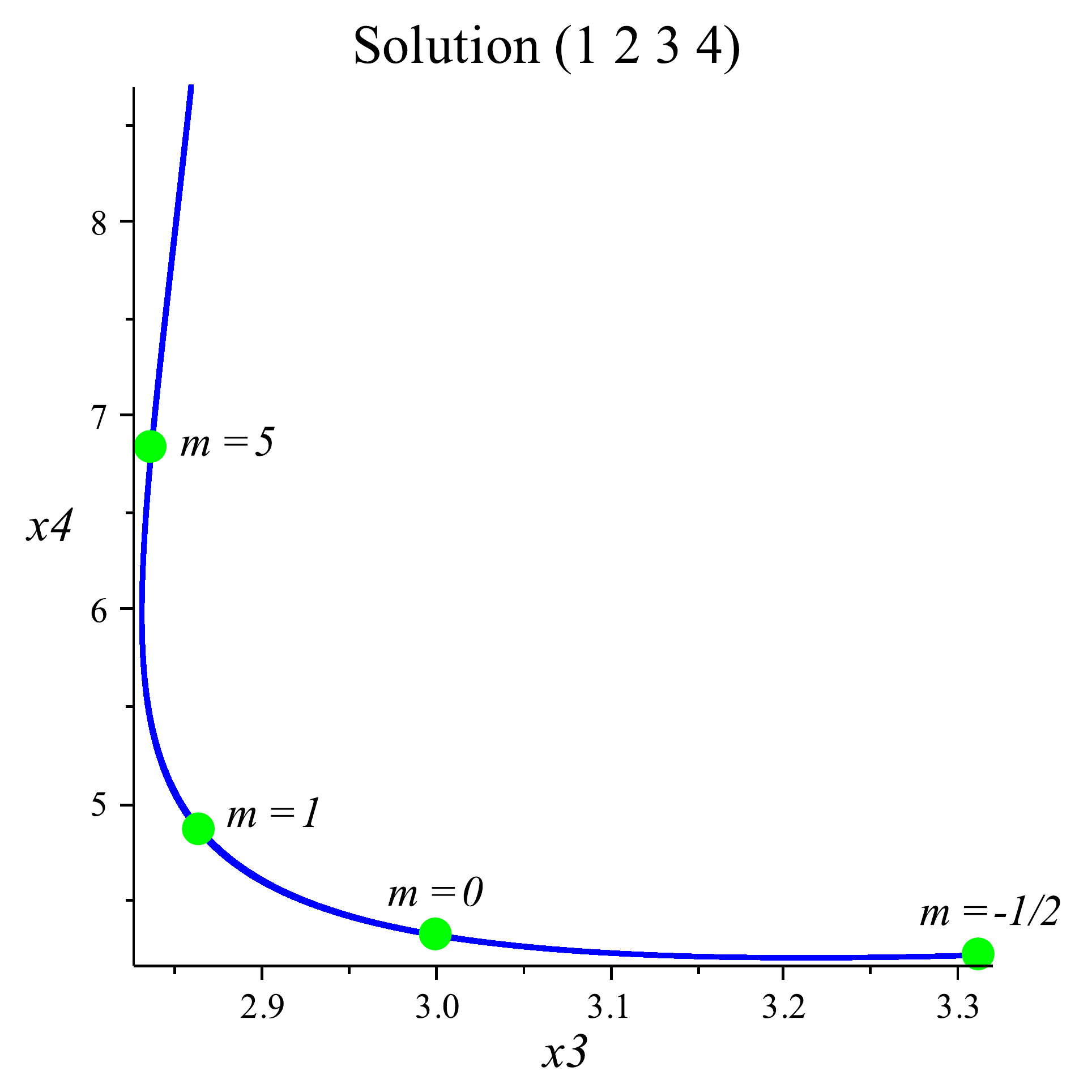}
\hspace{0.1in}
\includegraphics[width=220bp]{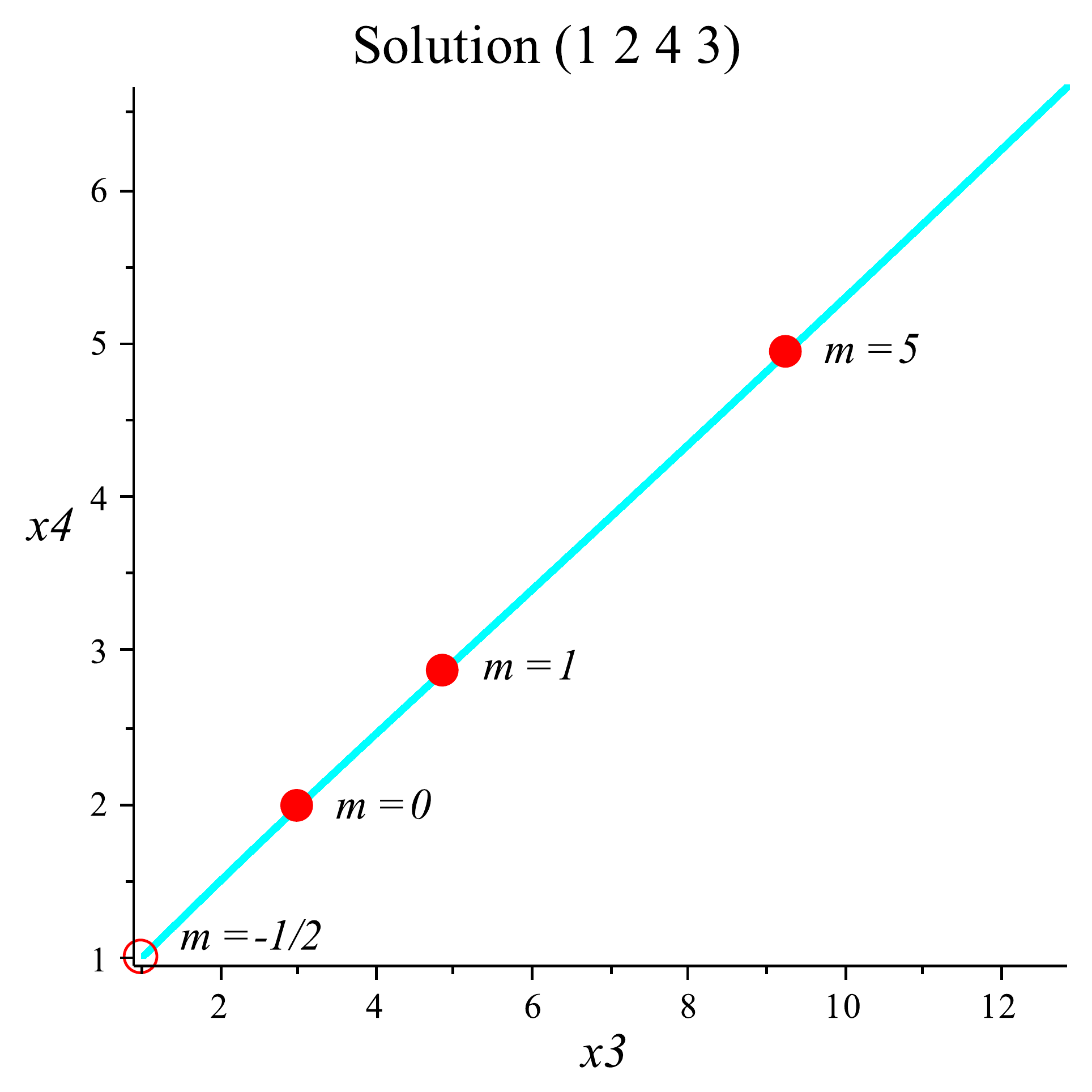}
\end{center}
\vspace{-0.2in}
\caption{A plot of $x_4$ versus $x_3$ for the collinear relative equilibrium with ordering $\ord{1}{2}{3}{4}$ (left) and $\ord{1}{2}{4}{3}$ (right) as the parameter $m$ varies from $-1/2$ to $10$.
Note the continuous transition as $m$ flips sign.}
\label{fig:sols}
\end{figure}

Some care must be taken to continue the solution with ordering $\ord{1}{2}{3}{4}$ to $m < 0$ because $\rho_i = 0 \; \forall i$ at $m = 0$ and $P$ has repeated roots.  
However, as explained in Example~\ref{ex:zero}, there are 12 different solutions for the case $m = 0$, one for each possible ordering.
A similar calculation to the one outlined in Example~\ref{ex:equal} shows that for the case $m = -1/4$, there are also 12 solutions, one for each ordering.
In this case, we take $\rho$ to be the smaller (in absolute value) {\em negative} root of $P_1$ in order to obtain the solution with $1 < x_3 < x_4$.
As explained above, this solution can be continued throughout the interval $-1/2 < m < 0$ because $x_3 = x_2 = 1$ and $x_3 = x_4$ are impossible.  
Thus, the solution with ordering $\ord{1}{2}{3}{4}$ varies continuously as $m$ decreases through $0$, as the corresponding choice for~$\rho$ transitions
from the smallest positive root of $P_1$ to the smallest negative root of $P_1$, passing through $\rho = 0$ at $m = 0$.  
Applying Theorem~\ref{ThmInv}, we have shown that the six orderings from Group~I are realized for any $m > -1/2$.
A plot of the solution curve in the $x_3x_4$-plane that corresponds to the ordering $\ord{1}{2}{3}{4}$ for $-1/2 \leq m \leq 10$ is shown to the left in Figure~\ref{fig:sols}.

The argument for the ordering $\ord{1}{2}{4}{3}$ and its five cousins is similar.   This time we follow the larger (in absolute value) negative root of $P_1$ for $m > 0$ because
$\rho = -\sqrt{3} - \sqrt{2}$ corresponds to the solution with $1 < x_4 < x_3$ when $m = 1$.  For $-1/2 < m < 0$, we follow the larger positive root of~$P_1$, 
making a continuous transition through $\rho = 0$ at $m = 0$.
We know the solution satisfying $1 < x_4 < x_3$ persists for all $m > -1/2$ because the only possible collisions are at $x_4 = x_3$ 
and $x_4 = 1$.  The first of these was eliminated by the Gr\"{o}bner basis calculation mentioned above, 
while the second is impossible because $f_5 = u(x_4 - 1) - 1$ was included in the original calculation of $GB$.  
By Theorem~\ref{ThmInv}, it follows that all six orderings from Group~II are realized.
A plot of the solution curve in the $x_3x_4$-plane that corresponds to the ordering $\ord{1}{2}{4}{3}$ for $-1/2 < m \leq 10$ is shown on the right in Figure~\ref{fig:sols}.
Although the curve appears to be linear, it is not.  This completes the proof of item {\bf (i)}.

\vs

{\bf (ii)}  A bifurcation occurs at $m = -1/2$ as the quartic $P_1$ becomes a quadratic with roots at $\pm \sqrt{3/7}$.  Here, $P$ reduces to a tenth degree polynomial with
repeated roots at $-1$ and $1$, each with multiplicity two.  These correspond to collisions between $x_3$ and $x_1$, or $x_3$ and $x_2$, respectively.
The remaining six roots of $P$ give six relative equilibria with the orderings from Group~I.  This can be shown rigorously by calculating the largest root of $P_2$ when
$\rho = -\sqrt{3/7}$ and then finding $x_4$ from $Q$.  We find that $1 < x_3 < x_4$, so this solution has ordering $\ord{1}{2}{3}{4}$.  
Theorem~\ref{ThmInv} then yields the remaining five solutions from Group~I.

For $-1 < m < -1/2$, Lemmas \ref{Lemma:QuarticRho} and~\ref{Lemma:Cubic} combine to show that $P$ has six distinct real roots and six complex roots.
By Lemma~\ref{lemma:extend}, the six real roots can be extended to a full solution of system~(\ref{eq:orgSys}).  As with case {\bf (i)}, rigorously justifying that the six
solutions belong to the orderings from Group~I involves picking a sample test case (we choose $m = -3/4$) and showing that the solution with ordering
$\ord{1}{2}{3}{4}$ persists for all $m$ in the open interval $(-1, -1/2)$.  Here we follow the negative root of $P_1$.
The argument is similar to that used in case {\bf (i)}.

\vs

{\bf (iii)}  If $m \leq -1$, then all of the roots of $P_1$ are complex.  Applying Lemma~\ref{Lemma:Cubic} shows that $P_2$ has no real roots, so there are no solutions.
\enpf

\begin{remark}
For $m > 0$, the existence of a unique collinear relative equilibrium for each possible ordering is a consequence of a well-known result from the Newtonian $n$-body problem
due to Moulton~\cite{moulton}.  For any choice of positive masses, there are exactly $n!/2$ collinear relative equilibria, one for each possible ordering
(see Section~2.1.5 of~\cite{meyer} or Section 2.9 of~\cite{rick-book}).  The result generalizes to the vortex setting as long as the circulations are
positive~\cite{palmore, oneil}.  For the case $-1 < m < 0$, a result due to O'Neil implies that there are at least six solutions (Theorem~6.2.1 in~\cite{oneil}).
\end{remark}

\subsection{Bifurcations}
\label{subsec:bifur}

We now discuss the bifurcations at $m = -1/2$ and $m = -1$ in greater detail, focusing on the behavior of those solutions which disappear after the bifurcation.     
As $m$ decreases toward $-1/2$, four of the solutions with orderings from Group II head toward
triple collision, while the remaining two orderings have the third vortex escaping to $\pm \infty$.  To see this, note that two of the roots of $P_1$ are heading off to $\pm \infty$ as
$m \rightarrow -1/2^{+}$.  Focusing on the ordering $\ord{1}{2}{4}{3}$, we track the solution corresponding to $\rho \rightarrow + \infty$.  Since $\rho$ is positive, we have $1 <  x_4 < x_3 < 3$
for this particular solution.  Moreover, since
$$
x_3^3 + \rho x_3^2  - 9x_3 - \rho = 0  \quad \Longrightarrow \quad  x_3^2 - 1 =  \frac{x_3(9 - x_3^2)}{\rho},
$$
we see that $x_3 \rightarrow1^{+}$ as $\rho \rightarrow \infty$.  By the Squeeze Theorem, we also have that $x_4 \rightarrow 1^{+}$, and thus the limiting
configuration has a triple collision between vortices 2, 3, and~4.  To track the other five solutions from Group~II, we use the formulas in Table~\ref{table:orders} and
take limits as $a \rightarrow 1^{+}$ and $b \rightarrow 1^{+}$.  The solution with ordering $\ord{1}{3}{4}{2}$ also limits on triple collision between vortices 2, 3, and~4.  
Orderings $\ord{3}{4}{1}{2}$ and $\ord{1}{4}{3}{2}$ limit on
triple collision between vortices 1, 3, and~4.  The solutions with orderings $\ord{1}{4}{2}{3}$ and $\ord{3}{1}{4}{2}$ have $x_3 \rightarrow \infty$ and
$x_3 \rightarrow -\infty$, respectively, but the $x_4$-coordinate takes the form $0/0$.

To determine the fate of the fourth vortex for the orderings $\ord{1}{4}{2}{3}$ and $\ord{3}{1}{4}{2}$,
we first compute an asymptotic expansion for $x_3$ and $x_4$ corresponding to the ordering $\ord{1}{2}{4}{3}$.
Introduce the small parameter $\epsilon$ by setting $m = -1/2 + \epsilon^2$, and let $\kappa = 1/\rho$ be a new variable.  
As $\epsilon \rightarrow 0$, we have $\kappa \rightarrow 0$, while
$x_3$ and $x_4$ each approach~1.  Rewriting $P_1 = 0$ and $P_2 = 0$ with~$\kappa$, we can expand $\kappa$ and $x_3$ in powers of~$\epsilon$.  
This, in turn, leads to an expansion for~$x_4$.  We find that
\begin{eqnarray}
\kappa  & = &  \frac{\sqrt{14}}{7} \, \epsilon  + \frac{289}{1029}  \sqrt{14} \, \epsilon^3 + \frac{24373}{43218}  \sqrt{14} \, \epsilon^5 + \mathcal{O}(\epsilon^7),  \nonumber \\[0.07in]
x_3 & = &  1 + \frac{4}{7} \sqrt{14} \, \epsilon  + \frac{8}{7} \, \epsilon^2 - \frac{20}{1029} \sqrt{14} \, \epsilon^3 - \frac{416}{1029} \, \epsilon^4 + \mathcal{O}(\epsilon^5),  \label{eq:x3Gr2As}  \\[0.07in]
x_4 & = &  1 + \frac{2}{7} \sqrt{14} \, \epsilon  + \frac{4}{7} \, \epsilon^2 - \frac{10}{1029} \sqrt{14} \, \epsilon^3 - \frac{320}{1029} \, \epsilon^4 + \mathcal{O}(\epsilon^5)   \label{eq:x4Gr2As}
\end{eqnarray}
are expansions for the solution with ordering $\ord{1}{2}{4}{3}$ near $m = -1/2$.  
Using these expansions, we have
$$
\frac{-2x_4 + x_3 + 1}{x_3 - 1} \; = \;  \frac{4}{147} \sqrt{14} \, \epsilon^3  + \mathcal{O}(\epsilon^4),
$$
which implies that $x_4 \rightarrow 0$ as $m \rightarrow -1/2^{+}$ for both solutions with orderings $\ord{1}{4}{2}{3}$ and $\ord{3}{1}{4}{2}$.
Notice that $x_4 - 1 \approx (1/2)(x_3 - 1)$, an observation which supports the nearly linear relationship shown in the right-hand graph of Figure~\ref{fig:sols}.
Finally, the expansions given above are perfectly valid for $\epsilon < 0$ as well.  In this case, they correspond to the solution with ordering $\ord{1}{3}{4}{2}$
near $m = -1/2$.

As $m$ decreases through $-1$, the solutions with orderings from Group I vanish; however, in this case, four of the limiting configurations end with a pair of
binary collisions.  As $m \rightarrow -1^{+}$, the remaining real roots of $P_1$ are heading off to $\pm \infty$.  
Focusing on the ordering $\ord{1}{2}{3}{4}$, we track the solution corresponding to $\rho \rightarrow - \infty$.  By Lemma~\ref{Lemma:Cubic} we have
$3 < x_3 < x_4$ for this solution.  Since $P_2(-\rho) = 8 \rho < 0$ and the leading coefficient of $P_2$ is positive, we see that $P_2$ has a root larger
than $-\rho$.  In fact, for $m$ close to $-1$, $x_3 = -\rho$ is an excellent approximation to this root.  Thus, for the ordering $\ord{1}{2}{3}{4}$, 
$x_3 \rightarrow \infty$, which implies $x_4 \rightarrow \infty$ as well.  
A similar fate occurs for the solution with ordering $\ord{4}{3}{1}{2}$, except that here, $x_3$ and $x_4$ approach $- \infty$ as $m \rightarrow -1^{+}$.

\begin{figure}[tb]
\begin{center}
\includegraphics[width=250bp]{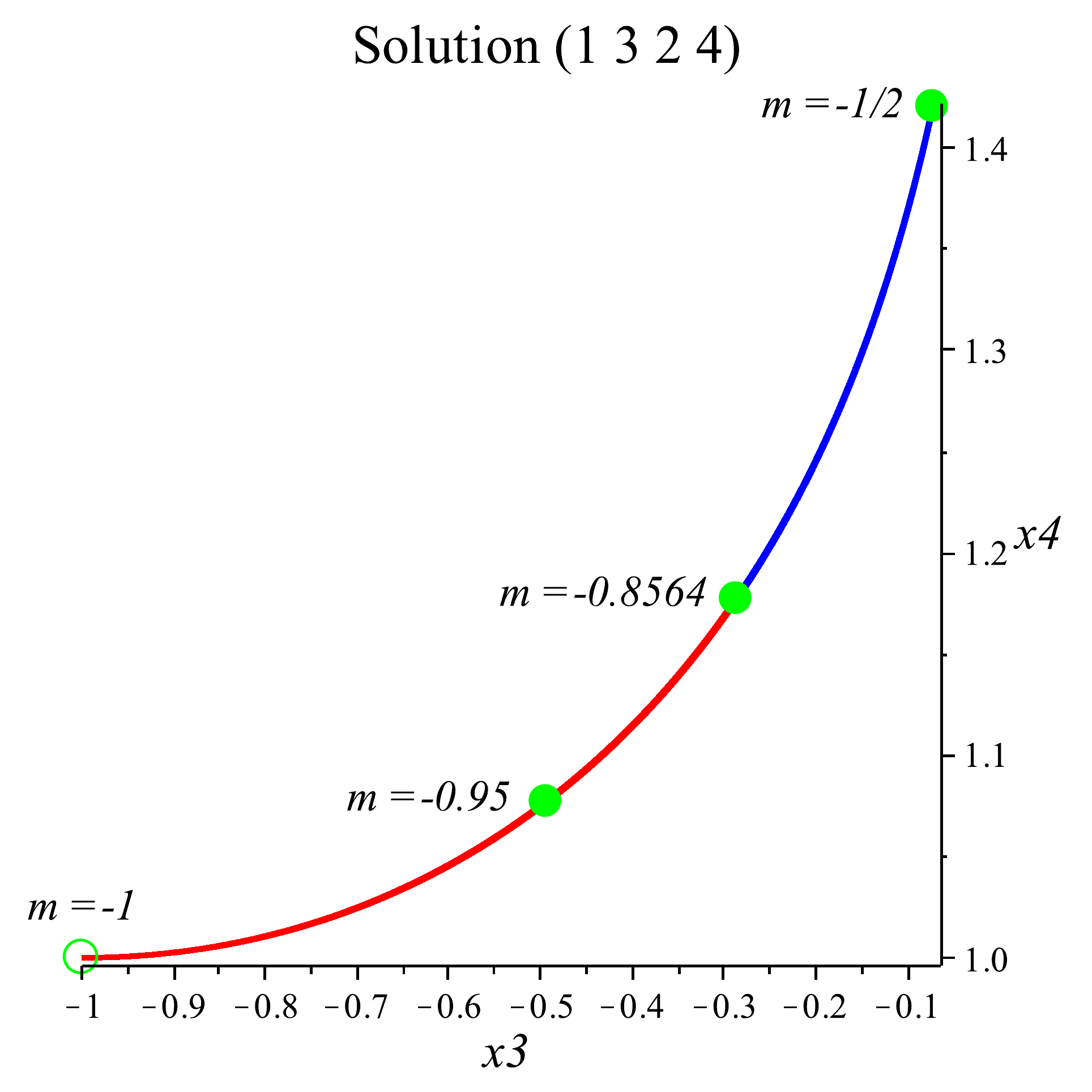}
\end{center}
\vspace{-0.2in}
\caption{A plot of $x_4$ versus $x_3$ for the collinear relative equilibrium with ordering $\ord{1}{3}{2}{4}$ as $m$ ranges from $-1/2$ to $-1$.  The red portion of the 
curve corresponds to linearly stable configurations, while the blue portion represents unstable solutions.  As $m \rightarrow -1^{+}$, vortex three approaches vortex one and
vortex four collides with vortex two.}
\label{fig:bifur}
\end{figure}

For the remaining four orderings in Group~I, the formulas in Table~\ref{table:orders} show that the $x_3$-coordinate is approaching $1$ or $-1$, while the $x_4$-coordinate
is an indeterminate form.  Substituting $m = -1$ and $x_3 = 1$ into $Q$ quickly yields $x_4 = -1$, while inserting $m = -1$ and $x_3 = -1$ into $Q$ gives
$x_4 = 1$.  Thus, as $m \rightarrow -1^{+}$, the solutions with orderings $\ord{4}{1}{2}{3}$ and $\ord{4}{1}{3}{2}$ have vortex four approaching vortex one,
and vortex three approaching vortex two.   For the orderings $\ord{1}{3}{2}{4}$ and $\ord{3}{1}{2}{4}$, the opposite collisions occur, 
with vortex three approaching vortex one, and vortex four approaching vortex two.   
A plot of the solution curve in the $x_3x_4$-plane that corresponds to the ordering $\ord{1}{3}{2}{4}$ for $-1 < m \leq -1/2$ is shown in Figure~\ref{fig:bifur}.
We will give an asymptotic expansion about $m = -1$ for this solution in Section~\ref{subsec:evalStructure}.

\section{Linear Stability of Solutions}

We now turn to the linear stability of the collinear relative equilibria found in Section~2, investigating the eigenvalues as the parameter~$m$ varies.
Through our analysis, we discover an important polynomial,
$$
\Psi \; = \;  64m^6 + 320m^5 +96m^4 - 220m^3 +505m^2 + 522m + 9, 
$$
whose roots include two new bifurcation values.  
Using Sturm's Theorem, $\Psi$ has four real roots, all of which are negative, and precisely two real roots between $-1$ and $0$.  
The root closest to $-1$ is  $m^{\ast} \approx -0.8564136$.  Note that $-6/7$ is a fairly good approximation to this root;  it is the second convergent
in the continued fraction expansion for $m^{\ast}$.   The root $m_c \approx -0.0175413$ will also be significant.
We will prove that the collinear relative equilibria in Group~I are linearly stable for $-1 < m < m^{\ast}$.   
For all other $m$-values, the relative equilibria in both groups are unstable.

\subsection{Background and a useful lemma}

We first review some key definitions and properties concerning the linear stability of a relative equilibrium~$z_0$ in the planar
$n$-vortex problem.    We follow the approach and setup described in~\cite{g:stability}.
The natural setting for determining the stability of $z_0$ is to change to rotating coordinates and treat 
$z_0$ as a rest point of the corresponding flow.   We will assume that $z_0$ has been translated so that its center of vorticity~$c$
is located at the origin.

Denote $M = {\rm diag}\{ \Gamma_1, \Gamma_1, \ldots, \Gamma_n, \Gamma_n \} $ as
the $2n \times 2n$ matrix of circulations, and let $K$ be the $2n \times 2n$ block diagonal matrix
containing 
$
J=
\begin{bmatrix}
0 & 1\\
-1 & 0
\end{bmatrix}
$
on the diagonal.  
The matrix that determines the linear stability of a relative equilibrium $z_0$ with angular velocity~$\omega$ is given by
$$
B \; = \;  K( M^{-1} D^2 H(z_0) + \omega I ) ,
$$
where $D^2H(z_0)$ is the Hessian of the Hamiltonian evaluated at $z_0$ and
$I$ is the $2n \times 2n$ identity matrix.
Since we are working with a Hamiltonian system, the eigenvalues of $B$ come in pairs $\pm \lambda$.  For a solution to be linearly stable,
the eigenvalues must lie on the imaginary axis.

One important property of the Hessian is that it anti-commutes with $K$, that is,
\begin{equation}
D^2H(z) K \; = \;  -K D^2H(z) .
\label{eq:commute}
\end{equation}
From this, it is straight-forward to see that
the characteristic polynomial of $M^{-1} D^2 H(z_0)$ is even.  In addition, if $v$ is an eigenvector 
of $M^{-1} D^2 H(z_0)$ with eigenvalue $\mu$, then $Kv$ is also an eigenvector with eigenvalue $- \mu$
(see Lemma~2.4 in~\cite{g:stability}).  This fact cuts the dimension of the problem in half.

In order to compare the collinear relative equilibria within a particular group, it is easier to work with the {\em scaled stability matrix}
\begin{equation}
\omega^{-1} B \; = \;  K( \omega^{-1} M^{-1} D^2H(z_0) + I) .  
\label{matrix-scaled}
\end{equation}
This scaling has no effect on the stability of $z_0$ because the characteristic polynomial of $B$ is even.
We will refer to the eigenvalues of $\omega^{-1} B$ as {\em normalized eigenvalues}.
The following lemma explains how to compute
the normalized eigenvalues from the eigenvalues of $\omega^{-1} M^{-1} D^2H(z_0)$.

\begin{lemma}
Let $p(\lambda)$ denote the characteristic polynomial of the scaled stability matrix $\omega^{-1} B$.
\begin{itemize}
\item[{\bf (i)}]
Suppose that $v$ is a real eigenvector of  $\omega^{-1} M^{-1} D^2H(z_0)$ with eigenvalue $\mu$.  Then $\{ v, Kv \}$ is a real invariant subspace of $ \omega^{-1} B$
and the restriction of $\omega^{-1} B$ to $\{ v, Kv \}$ is
\begin{equation}
\begin{bmatrix}
0 &\mu - 1 \\
\mu +1 & 0 \\
\end{bmatrix}.
\label{eq:mat1}
\end{equation}
Consequently, $p(\lambda)$ has a quadratic factor of the form $\lambda^2 + 1- \mu^2$.

\item[{\bf (ii)}]
Suppose that $v = v_1 + i v_2$ is a complex eigenvector of  $\omega^{-1} M^{-1} D^2H(z_0)$ 
with complex eigenvalue $\mu = \alpha + i \beta$.  Then $\{ v_1,v_2, Kv_1, Kv_2 \}$ is a real invariant subspace of $\omega^{-1} B$ 
and the restriction of $\omega^{-1} B$ to this space is
\begin{equation}
\begin{bmatrix}
0 & 0 & \alpha - 1 & \beta \\
0 & 0 & -\beta &  \alpha - 1 \\
\alpha + 1 & \beta & 0 & 0 \\
-\beta &\alpha + 1 & 0 & 0 
\end{bmatrix}.
\label{eq:mat2}
\end{equation}
Consequently, $p(\lambda)$ has a quartic factor of the form 
$(\lambda^2 + 1 - \mu^2)(\lambda^2 + 1 - \overline{\mu}^{\, 2})$, where $\overline{\mu} = \alpha - i \beta$.
\end{itemize}

\label{lemma:invsubspace}
\end{lemma}

\pf
{\bf (i)} Since $\omega^{-1} M^{-1} D^2H(z_0) v = \mu v$, we have  
$\omega^{-1} M^{-1} D^2H(z_0) K v = - \mu K v$ by equation~(\ref{eq:commute}).
This implies that
$\omega^{-1} B v = (\mu + 1) K v$ and $\omega^{-1} B (Kv) =  (\mu -1) v$, verifying matrix~(\ref{eq:mat1}).  The characteristic polynomial
of matrix~(\ref{eq:mat1}) is $\lambda^2 + 1 - \mu^2$ and therefore, this quadratic is a factor of $p(\lambda)$.

{\bf (ii)}  If $v_1 + i v_2$ is a complex eigenvector with eigenvalue $\alpha + i \beta$, then we have
$\omega^{-1} M^{-1} D^2H(z_0) v_1 = \alpha v_1 - \beta v_2$ and
$\omega^{-1} M^{-1} D^2H(z_0) v_2 = \alpha v_2 + \beta v_1$.
Using equation~(\ref{eq:commute}), this implies that
$\omega^{-1} B v_1 = (\alpha+ 1) K v_1 - \beta K v_2$, $\omega^{-1} B v_2 =  \beta K v_1 + (\alpha + 1) K v_2$,
$\omega^{-1} B (Kv_1) =  (\alpha - 1) v_1 - \beta v_2$, and $\omega^{-1} B (Kv_2) =  \beta v_1 + (\alpha - 1) v_2$,
which confirms matrix~(\ref{eq:mat2}).  The characteristic polynomial
of matrix~(\ref{eq:mat2}) is $(\lambda^2 + 1 - \mu^2)(\lambda^2 + 1 - \overline{\mu}^{\, 2})$
and hence, this quartic is a factor of $p(\lambda)$.
\enpf

Due to the conserved quantities of the $n$-vortex problem, any relative equilibrium will have the
four normalized eigenvalues  $0, 0, \pm i$.  We call these eigenvalues {\em trivial}. 
The eigenvalues $\pm i$ arise from the center of vorticity integral and can be derived
by noting that the vector $[1, 0, 1, 0, \ldots, 1, 0]$ is in the kernel of 
$\omega^{-1} M^{-1} D^2H(z_0)$.  The two zero eigenvalues appear because
relative equilibria are not isolated rest points.  For any relative equilibrium $z_0$, the
vector $K z_0$ is in the kernel of $\omega^{-1} B$.   This vector is tangent to the periodic orbit 
determined by $z_0$ at $t=0$.  This follows from the identity
\begin{equation}
\omega^{-1} M^{-1} D^2H(z_0) \, z_0 \; = \;  z_0,
\label{eq:evec-RE}
\end{equation}
and part~{\bf (i)} of Lemma~\ref{lemma:invsubspace}.  
Thus, in the full phase space, a relative equilibrium is always degenerate.

One method for dealing with the issues that arise from the symmetries of the problem is to work
in a reduced phase space (e.g., quotienting out the rotational symmetry).   However, it is typically easier to make
the computations in $\mathbb{R}^{2n}$ and then define linear stability by restricting to the appropriate subspace.
This is the approach we follow here.

Let $V = \mbox{span} \{z_0, Kz_0\}$ and denote $V^\perp$ as the $M$-orthogonal complement of $V$, that is,
$$
V^\perp \; = \;  \{ w \in  \mathbb{R}^{2n} : w^T M v = 0 \; \; \forall v \in V \} .
$$
The invariant subspace $V$ accounts for the two zero eigenvalues;
the vector space $V^\perp$ has dimension $2n - 2$ and is invariant under $\omega^{-1} B$.  
We also have that $V \cap V^\perp = \{0\}$ provided $L \neq 0$.  This motivates the following definition
for linear stability.

\begin{define}
A relative equilibrium $z_0$ always has the four trivial normalized eigenvalues $0, 0, \pm i$.  We call $z_0$ {\em nondegenerate}
if the remaining $2n - 4$ eigenvalues are nonzero.  A nondegenerate relative equilibrium  is {\em spectrally stable} if the nontrivial 
eigenvalues lie on the imaginary axis, and {\em linearly stable} if, in addition, the restriction of the scaled stability matrix $\omega^{-1} B$ to $V^\perp$ 
has a block-diagonal Jordan form with blocks
$
\begin{bmatrix}
0 &  \beta_i \\  
-\beta_i & 0 \\
\end{bmatrix}.
$
\end{define}

As noted in~\cite{g:stability}, if $\Gamma_i > 0 \, \forall i$, then $\omega^{-1} M^{-1} D^2H(z_0)$ is symmetric with respect to an $M$-orthonormal
basis and has a full set of linearly independent real eigenvectors.  Consequently, part~{\bf (i)} of
Lemma~\ref{lemma:invsubspace} applies repeatedly and the characteristic polynomial factors as
$$
p(\lambda) \; = \;  \lambda^2 (\lambda^2 + 1) \prod_{j = 1}^{n-2} (\lambda^2 + 1 - \mu_j^2),
$$
where the $\mu_j$ are the nontrivial eigenvalues of $\omega^{-1} M^{-1} D^2H(z_0)$.
It follows that the relative equilibrium is linearly stable if and only if $|\mu_j| < 1 \; \forall j$.

If the circulations $\Gamma_i$ are of mixed sign, then $\omega^{-1} M^{-1} D^2H(z_0)$ 
may have complex eigenvalues, leading to quartic factors of the characteristic polynomial, as explained in part~{\bf (ii)}
of Lemma~\ref{lemma:invsubspace}.  This is the case for certain values of $m < 0$ in our problem.  Note that when $\mu_j \in \mathbb{C} - \mathbb{R}$,
the corresponding eigenvalues of the relative equilibrium form a complex quartuplet $\pm \alpha' \pm i \beta'$.  This implies instability unless
$\alpha' = 0$, which only occurs when Re$(\mu_j) = 0$.

\subsection{Finding the nontrivial eigenvalues of a collinear relative equilbrium}

We now focus on the stability of a collinear relative equilibrium $z_0 = (z_1, z_2, \ldots, z_n) \in \mathbb{R}^{2n}$,
where $z_i = (x_i, 0) \; \forall i$.  Rearranging the coordinates from $(x_1, y_1, x_2, y_2, \ldots, x_n, y_n)$ to 
$(x_1, x_2, \ldots, x_n, y_1, y_2, \ldots, y_n)$, we find that $\omega^{-1} M^{-1} D^2H(z_0)$ takes the special form
$$
\omega^{-1} M^{-1} D^2H(z_0) \; = \;  
\begin{bmatrix}
A & 0  \\
0 &  -A
\end{bmatrix} ,
$$
where $A$ is an $n \times n$ matrix with entries $A_{ij} = - \omega^{-1} \Gamma_j  r_{ij}^{-2}$ if $i \neq j$ and
$A_{ii} = - \sum_{j \neq i}  A_{ij}$.  This reduces our calculations from a $2n$-dimensional vector space to an $n$-dimensional
one.

For the case $n=4$, we have
$$
A = \omega^{-1}
\begin{bmatrix}
\dsty{\frac{ \Gamma_2}{r_{12}^2} + \frac{ \Gamma_3}{r_{13}^2} + \frac{  \Gamma_4}{r_{14}^2} }
&  \dsty{ -\frac{ \Gamma_2}{r_{12}^2} } & \dsty{ -\frac{  \Gamma_3}{r_{13}^2}} & \dsty{ -\frac{  \Gamma_4}{r_{14}^2} }  \\[0.25in]
\dsty{-\frac{\Gamma_1}{r_{12}^2}}  &  \dsty{\frac{\Gamma_1}{r_{12}^2} + \frac{\Gamma_3}{r_{23}^2} + \frac{\Gamma_4}{r_{24}^2} }
& \dsty{ -\frac{\Gamma_3}{r_{23}^2}} &  \dsty{-\frac{\Gamma_4}{r_{24}^2}  } \\[0.25in]
\dsty{ -\frac{\Gamma_1}{r_{13}^2}}  & \dsty{ -\frac{\Gamma_2}{r_{23}^2}} & 
\dsty{ \frac{\Gamma_1}{r_{13}^2} + \frac{\Gamma_2}{r_{23}^2} + \frac{\Gamma_4}{r_{34}^2}}  &  \dsty{-\frac{\Gamma_4}{r_{34}^2} }  \\[0.25in]
\dsty{ -\frac{\Gamma_1}{r_{14}^2} } & \dsty{ -\frac{\Gamma_2}{r_{24}^2}} & \dsty{ -\frac{\Gamma_3}{r_{34}^2}}  &   
\dsty{ \frac{\Gamma_1}{r_{14}^2} + \frac{\Gamma_2}{r_{24}^2} + \frac{\Gamma_3}{r_{34}^2} }
\end{bmatrix} .
$$
Note that the vectors $s = [1, 1, 1, 1]^{\mbox{\small{T}}}$ and $x = [x_1, x_2, x_3, x_4]^{\mbox{\small{T}}}$ are eigenvectors of $A$ with eigenvalues 0 and~1, respectively.
These are the two trivial eigenvalues of $A$ arising from the center of vorticity integral and the rotational symmetry.  The remaining two eigenvalues of $A$
determine the linear stability of~$z_0$.  Specifically, applying both parts of Lemma~\ref{lemma:invsubspace},
the normalized nontrivial eigenvalues of $z_0$ are given by
\begin{equation}
\lambda_1 \; = \;  \pm  \sqrt{\mu_1^2 - 1}  \quad \mbox{ and} \quad
\lambda_2 \; = \;  \pm  \sqrt{\mu_2^2 - 1} \,  ,
\label{eq:nontrivEvals}
\end{equation}
where $\mu_1$ and $\mu_2$ are the nontrivial eigenvalues of~$A$.

Let $W = \mbox{span}\{s, x\}$ and let $W^\perp$ denote the $M$-orthogonal complement of $W$
where $M = \break  {\rm diag}\{ \Gamma_1, \Gamma_2, \Gamma_3, \Gamma_4 \} $.  The subspace $W^\perp$ is invariant under $A$.
To find $\mu_1$ and $\mu_2$, we compute the restriction of $A$ to~$W^\perp$.  
It is straight-forward to check that the two vectors
\begin{eqnarray*}
w_1 & = &  [\Gamma_2 \Gamma_3 (x_3 - x_2),  \; \Gamma_1 \Gamma_3 (x_1 - x_3),  \; \Gamma_1 \Gamma_2 (x_2 - x_1),  \; 0]^{\mbox{\small{T}}},  \\[0.07in]
w_2 & = &  [\Gamma_2 \Gamma_4 (x_4 - x_2),  \; \Gamma_1 \Gamma_4 (x_1 - x_4),  \;  0, \;  \Gamma_1 \Gamma_2 (x_2 - x_1)]^{\mbox{\small{T}}}
\end{eqnarray*}
form a basis for $W^\perp$, as $w_i^{\mbox{\small{T}}} M s = 0$ and $w_i^{\mbox{\small{T}}} M x = 0$ for each~$i$.  However,
it is not an $M$-orthogonal basis since $w_1^{\mbox{\small{T}}} M w_2 \neq 0$.

Let $C$ denote the restriction of $A$ to~$W^\perp$ and write
$$
C \; = \;  
\begin{bmatrix}
C_{11} & C_{12}  \\
C_{21} &  C_{22}
\end{bmatrix} .
$$
To find the entries of $C$, note that $C_{11} w_1 + C_{21} w_2 = [\ast, \ast, C_{11} \Gamma_1 \Gamma_2 (x_2 - x_1), C_{21} \Gamma_1 \Gamma_2 (x_2 - x_1)]^{\mbox{\small{T}}}$.
Thus, we see that $C_{11} \Gamma_1 \Gamma_2 (x_2 - x_1)$ and $C_{21} \Gamma_1 \Gamma_2 (x_2 - x_1)$ are equal to the third and fourth coordinates, respectively, 
of the vector $Aw_1$.  A similar fact applies for $C_{12}$ and $C_{22}$.  
After some computation, we find that
\begin{eqnarray*}
C_{11} & = & \omega^{-1} \left[ \frac{\Gamma_1 + \Gamma_3}{r_{13}^2} + \frac{\Gamma_2 + \Gamma_3}{r_{23}^2} + \frac{\Gamma_4}{r_{34}^2} + \frac{\Gamma_3}{(x_3 - x_1)(x_3 - x_2)} \right], \\[0.07in]
C_{22} & = &   \omega^{-1} \left[ \frac{\Gamma_1 + \Gamma_4}{r_{14}^2} + \frac{\Gamma_2 + \Gamma_4}{r_{24}^2} + \frac{\Gamma_3}{r_{34}^2} + \frac{\Gamma_4}{(x_4 - x_1)(x_4 - x_2)} \right], \\[0.07in]
C_{21} & = &  - \omega^{-1} \frac{ \Gamma_3}{x_2 - x_1} \left[  \frac{x_3 - x_2}{r_{14}^2} + \frac{x_1 - x_3}{r_{24}^2} + \frac{x_2 - x_1}{r_{34}^2} \right], \\[0.07in]
C_{12} & = &  - \omega^{-1} \frac{ \Gamma_4}{x_2 - x_1} \left[  \frac{x_4 - x_2}{r_{13}^2} + \frac{x_1 - x_4}{r_{23}^2} + \frac{x_2 - x_1}{r_{34}^2} \right]. 
\end{eqnarray*}

The nontrivial eigenvalues of $A$, $\mu_1$ and $\mu_2$, are equivalent to the eigenvalues of $C$.  
They are easily expressed in terms of the trace and determinant of $C$.
The quantity~$\delta$, defined by 
$$
\delta  =  \frac{ (\Gamma_1 + \Gamma_2)(\Gamma_3 + \Gamma_4) }{r_{12}^2 r_{34}^2} +  \frac{ (\Gamma_1 + \Gamma_3)(\Gamma_2 + \Gamma_4) }{r_{13}^2 r_{24}^2} 
+ \frac{ (\Gamma_1 + \Gamma_4)(\Gamma_2 + \Gamma_3) }{r_{14}^2 r_{23}^2}  + \sum_{i=1}^4 \sum_{\substack{ j < k \\ j, k \neq i}}^4 \frac{ \Gamma_i (\Gamma_i + \Gamma_j + \Gamma_k) }{r_{ij}^2 r_{ik}^2},
$$
is important in the computation of the determinant of~$C$.

\begin{lemma}
Let $T$ and $D$ denote the trace and determinant, respectively, of $C$.  We have
\begin{eqnarray}
T  & = &  \omega^{-1}  \sum_{i < j}^4  \frac{\Gamma_i + \Gamma_j}{r_{ij}^2} \; - \;  1 ,  \label{eq:trace} \\
D & = &  -T + \omega^{-2} \delta.  \label{eq:det}
\end{eqnarray}
\end{lemma}

\pf  
By definition, if $(x_1, x_2, x_3, x_4)$ are the coordinates of a collinear relative equilibrium with center of vorticity at the origin, 
then equation~(\ref{eq:rel-equ}) implies that
\begin{eqnarray}
\frac{\Gamma_2}{x_2 - x_1} + \frac{\Gamma_3}{x_3 - x_1} + \frac{\Gamma_4}{x_4 - x_1} + \omega x_1 & = & 0,  \label{eq:cc1} \\[0.07in]
\frac{\Gamma_1}{x_1 - x_2} + \frac{\Gamma_3}{x_3 - x_2} + \frac{\Gamma_4}{x_4 - x_2} + \omega x_2 & = & 0.  \label{eq:cc2} 
\end{eqnarray}
Subtracting equation~(\ref{eq:cc2}) from equation~(\ref{eq:cc1}) and multiplying through by $1/(x_2 - x_1)$ gives
$$
\omega \; = \;    \frac{\Gamma_1 + \Gamma_2}{r_{12}^2}  - \frac{\Gamma_3}{(x_3 - x_1)(x_3 - x_2)} - \frac{\Gamma_4}{(x_4 - x_1)(x_4 - x_2)} .
$$
It then follows that
$$
T \; = \;  C_{11} + C_{22} \; = \;   \omega^{-1}  \sum_{i < j}^4  \frac{\Gamma_i + \Gamma_j}{r_{ij}^2} \; - \;  1.
$$

Alternatively, recall that the trace of a matrix is equal to the sum of its eigenvalues.  Applying this fact to both matrices $A$ and $C$ gives
$$
 1 + T \; = \;  0 + 1 + \mu_1 + \mu_2   \; = \;  \mbox{Tr}(A) \; = \;  \omega^{-1}  \sum_{i < j}^4  \frac{\Gamma_i + \Gamma_j}{r_{ij}^2} ,
$$
which gives an alternative proof of formula~(\ref{eq:trace}).

Computing the determinant $D$ from the entries of $C$ gives a messy expression.  A more useful formula can be obtained by utilizing the fact that
the sum of the product of all pairs of eigenvalues of~$A$ is equal to half the quantity $(\mbox{Tr}(A))^2 - \mbox{Tr}(A^2)$.  This yields
$$
0 \cdot 1 + 0 \cdot \mu_1 + 0 \cdot \mu_2 + 1 \cdot \mu_1 + 1 \cdot \mu_2 + \mu_1 \cdot \mu_2 \; = \;  \frac{1}{2} \left[ 
\omega^{-2} \left(  \sum_{i < j}^4  \frac{\Gamma_i + \Gamma_j}{r_{ij}^2} \right)^2 - \mbox{Tr}(A^2) \right],
$$
which, after some calculation, gives
$$
T + D \; = \;  \omega^{-2} \delta \, .
$$
\enpf

\begin{remark}
Recall that the angular velocity of a relative equilibrium is given by $\omega = L/I$, where $L = \sum_{i < j} \Gamma_i \Gamma_j$ and 
$I = (1/\Gamma) \sum_{i < j} \Gamma_i \Gamma_j r_{ij}^2$.
It follows from formulas (\ref{eq:trace}) and~(\ref{eq:det}), that $T$ and $D$ depend only on the circulations~$\Gamma_i$ and the mutual distances~$r_{ij}$.
Thus, as we would expect, the stability of a relative equilibrium is unaffected by translation, and we may retain our
original coordinates (e.g., $x_1 = -1, x_2 = 1$) when calculating $T$ and~$D$, rather than shifting the configuration so that
the center of vorticity is at the origin.
\end{remark}

\begin{theorem}
The nontrivial eigenvalues of $A$ are the roots of $\lambda^2 - T \lambda + D$, where $T$ and $D$ are given by
(\ref{eq:trace}) and~(\ref{eq:det}), respectively.   They are identical for any solution within a particular
group of orderings.  Sufficient conditions for linear stability of the relative equilibrium are 
\begin{center}
{\bf (i)} $-2 < T < 2$, \qquad  {\bf (ii)} $D < T^2/4$,  \qquad  {\bf (iii)}  $D > T - 1$, \quad and \quad {\bf (iv)} $D > -T - 1$.
\end{center}
\label{thm:TDinv}
\end{theorem}

\pf
As explained above, the nontrivial eigenvalues of $A$ are equivalent to the eigenvalues of~$C$, and the characteristic
polynomial of $C$ is $\lambda^2 - T \lambda + D$.  It is straight-forward to check that $T$ and $D$ are invariant under the
maps $S$ and $R$ defined in equations (\ref{eq:symmS}) and~(\ref{eq:symmR}), respectively.  This was also confirmed using Maple.
It follows that $T$ and $D$ are invariant under the group~$G$ and thus, the nontrivial normalized eigenvalues are identical
for all six solutions in a given group of orderings.  

\begin{figure}[bth]
\begin{center}
\includegraphics[width=330bp]{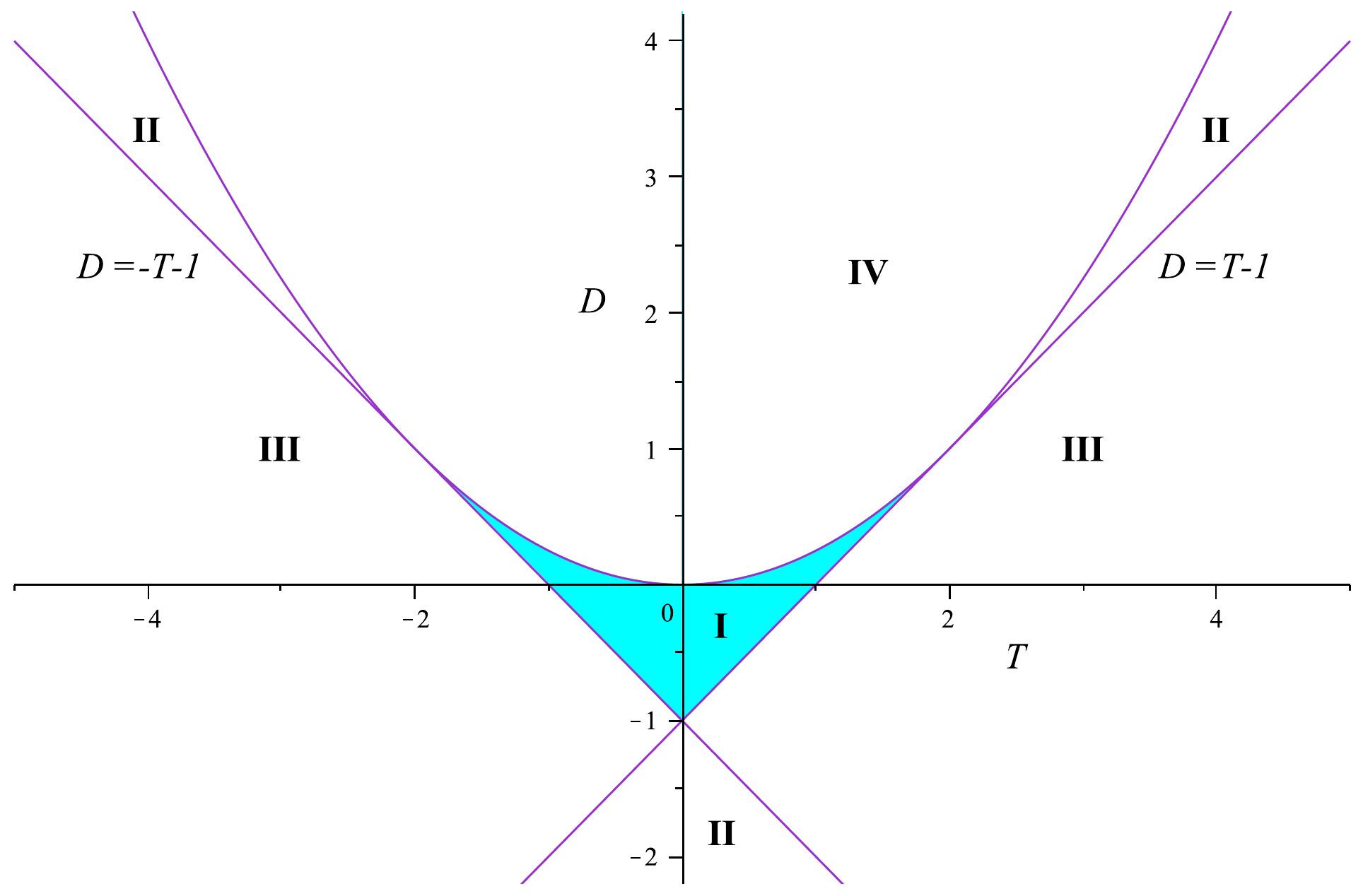}
\end{center}
\vspace{-0.2in}
\caption{The stability diagram in terms of the trace $T$ and determinant $D$ of $C$.  The light blue region 
corresponds to linearly stable solutions.  Zero eigenvalues occur
on the bifurcation lines $D = T - 1$ or $D = -T - 1$, while repeated eigenvalues arise on the parabola $D = T^2/4$.  
In region~II, there are two real pairs of eigenvalues $\pm \lambda_1, \pm \lambda_2$; 
in region~III, there is a real pair $\pm \lambda_1$ and a pure imaginary pair $\pm i \lambda_2$ of eigenvalues;
in region~IV, the eigenvalues form a complex quartuplet $\pm \alpha' \pm i \beta'$.}
\label{fig:TD-plane}
\end{figure}

If $\mu_i \in \mathbb{R}$, then formula~(\ref{eq:nontrivEvals}) shows that $|\mu_1| < 1$ and  $|\mu_2| < 1$ are both required for stability.  
This is guaranteed if conditions {\bf (i)} though~{\bf (iv)} are satisfied (see Figure~\ref{fig:TD-plane}).  Note that linear stability follows as well because 
condition~{\bf (ii)} insures that $\mu_1 \neq \mu_2$, so there are no repeated nontrivial eigenvalues. 
\enpf

\begin{remark}
In addition to the blue region in Figure~\ref{fig:TD-plane}, solutions are also linearly stable along the positive $D$-axis ($T = 0, D > 0$).  
In this case the eigenvalues of $C$ are pure imaginary and the nontrivial eigenvalues 
of the relative equilibrium are $\pm i \sqrt{1 + D} \, , \pm i \sqrt{1 + D} \,$.   Using matrix~(\ref{eq:mat2}), it is straight-forward to show that 
the Jordan form of $\omega^{-1} B$ has no off-diagonal blocks, and thus the solution is linearly stable.  
\end{remark}

\begin{example}
{\bf The case $m=1$}.

Recall from Example~\ref{ex:equal} that in the case of equal-strength vortices, the solution with ordering $\ord{1}{2}{3}{4}$ has positions $(x_1, x_2, x_3, x_4) = 
(-1, 1, -1 + \sqrt{2} + \sqrt{6}, 1 + \sqrt{2} + \sqrt{6})$.  This gives $T = 5$ and $D = 6$, so the nontrivial eigenvalues of $A$ are $\mu_1 = 2$ and $\mu_2 = 3$
and the relative equilibrium is unstable.  The same result holds for the ordering $\ord{1}{2}{4}{3}$.   By formula~(\ref{eq:nontrivEvals}),
the nontrivial normalized eigenvalues for either group are $\pm \sqrt{3}$ and $\pm 2\sqrt{2}$.  
\label{ExampleM1Stab}
\end{example}

\begin{example}
{\bf The case $m=0$}.

If $\Gamma_4 = 0$, we found in Example~\ref{ex:zero} that the solution with ordering $\ord{1}{2}{3}{4}$ has positions $(x_1, x_2, x_3, x_4) = 
(-1, 1, 3, 1 + (1/3)\sqrt{ 54 + 6 \sqrt{57} } \, )$.  From these values, we find that $\mu_1 = 2$ and $\mu_2 = (15 - \sqrt{57} \,)/4$, so
the relative equilibrium is unstable.  The nontrivial normalized eigenvalues are $\pm \sqrt{3} \approx \pm 1.732$ and
$\pm (1/4) \sqrt{ 266 - 30\sqrt{57}} \approx  \pm 1.571$.  A similar result holds for the solution with ordering $\ord{1}{2}{4}{3}$, except that
the eigenvalues are much further apart.  The nontrivial normalized eigenvalues for this solution are $\pm \sqrt{3}$ and 
$\pm (1/4) \sqrt{ 266 + 30\sqrt{57}} \approx \pm 5.548$ . 
\end{example}

\subsection{Stability, bifurcations, and eigenvalue structure}
\label{subsec:evalStructure}

We now apply Theorem~\ref{thm:TDinv} and formulas (\ref{eq:trace}) and~(\ref{eq:det}) to investigate the linear stability of our two families of relative equilibria 
as $m$~varies.  Recall from Theorem~\ref{thm:main} that there are two groups of solutions varying continuously in~$m$:
Group~I exists for all $m > -1$ and Group~II exists for all $m > -1/2$.  Moreover, for a fixed~$m$, the values 
of $x_3$ and $x_4$, and thereby the trace~$T$ and determinant~$D$, can be found analytically by working with the roots of a cubic equation.

We first compute an asymptotic expansion for the solution with ordering $\ord{1}{3}{2}{4}$ about $m = -1$.  Recall that as $m \rightarrow -1^{+}$,
the limiting configuration for this specific ordering contains a pair of binary collisions, that is, $x_3 \rightarrow -1$ and $x_4 \rightarrow 1$ (see Figure~\ref{fig:bifur}).
Following the same approach used in Section~\ref{subsec:bifur}, we let $m = -1 + \epsilon^2$, where $\epsilon$ is a small positive parameter, and
repeatedly solve the equations $P_1 = 0$, $P_2 = 0$, and $Q = 0$ to obtain the expansions
\begin{eqnarray}
x_3 & = &  -1 + 2\sqrt{2} \, \epsilon  -2 \epsilon^2 - 3 \sqrt{2} \, \epsilon^3 + 7 \epsilon^4 + \cdots - \frac{71179}{16} \sqrt{2} \, \epsilon^{11} + \frac{226695}{16} \epsilon^{12} + \mathcal{O}(\epsilon^{13}), \label{eq:x3As} \\[0.07in]
x_4 & = &  1 + 2 \epsilon^2  - \sqrt{2} \, \epsilon^3 - 5 \epsilon^4 + \frac{9}{2} \sqrt{2} \, \epsilon^5 + \cdots - \frac{24863}{16} \sqrt{2} \, \epsilon^{11} - \frac{192805}{16} \epsilon^{12} + \mathcal{O}(\epsilon^{13}).  \label{eq:x4As}
\end{eqnarray}
Substituting these expressions into the formula $\omega = L/I$ shows that $1/\omega = 4 \epsilon^2 + \mathcal{O}(\epsilon^3)$, which implies that the angular velocity becomes infinite as the
vortices approach collision.   Using formulas (\ref{eq:trace}) and~(\ref{eq:det}), we find the following expansions for $T$ and $D$:
\begin{eqnarray}
T & = & 1 + 6 \epsilon^2 - 12 \epsilon^4 + 60 \epsilon^6 - 426 \epsilon^8 +  \mathcal{O}(\epsilon^{10}),  \label{eq:TrAs}  \\[0.07in]
D & = &  6 \epsilon^2 - 12 \epsilon^4 + 78 \epsilon^6 - 588 \epsilon^8 + \frac{9453}{2} \epsilon^{10} + \mathcal{O}(\epsilon^{12}).   \label{eq:DtAs}
\end{eqnarray}
The fact that each series contains only even powers of $\epsilon$ is a consequence of the invariance described in Theorem~\ref{thm:TDinv}.
Choosing $\epsilon < 0$ in formulas (\ref{eq:x3As}) and~(\ref{eq:x4As}) is perfectly valid; in fact, it provides an expansion for the solution with 
ordering $\ord{3}{1}{2}{4}$, which is also a member of the Group~I orderings.  Since $T$ and $D$ are invariant within a specific group,
they must be even functions of the parameter~$\epsilon$.


Recall that $m^{\ast} \approx -0.8564136$ and $m_c \approx -0.0175413$ are important roots of the polynomial 
$$
\Psi(m) \; = \;  64m^6 + 320m^5 +96m^4 - 220m^3 +505m^2 + 522m + 9. 
$$

\begin{theorem}
The linear stability and nontrivial eigenvalue structure for the four-vortex collinear relative equilibria with circulations
$\Gamma_1 = \Gamma_2 = \Gamma_3 = 1$ and $\Gamma_4 = m$ are as follows:

\begin{itemize}
\item[{\bf (i)}]   The solutions from Group~I are linearly stable for $-1 < m < m^{\ast}$, spectrally stable at $m = m^{\ast}$, and unstable for $m > m^{\ast}$.
\item[{\bf (ii)}]   For $m > m_c$, the nontrivial eigenvalues for the Group~I solutions consist of two real pairs, and at $m = m_c$, these pairs merge to form a real pair with multiplicity two.  
As $m \rightarrow \infty$, the normalized Group~I eigenvalues approach $\pm 2 \sqrt{6}$ and $\pm 2 \sqrt{2}$.
For $m^{\ast} < m < m_c$, the nontrivial eigenvalues form a complex quartuplet $\pm \alpha \pm i \beta$.  As $m \rightarrow -1^{+}$, the nontrivial normalized eigenvalues approach
$0, 0, \pm i$.
\item[{\bf (iii)}]  The solutions from Group~II ($m > -1/2$) are always unstable with two real pairs of nontrivial normalized eigenvalues $\pm \lambda_1, \pm \lambda_2$.
As $m \rightarrow -1/2^{+}$, $\lambda_1 \rightarrow \infty$ and $\lambda_2 \rightarrow 2\sqrt{14}/5$.  
As $m \rightarrow \infty$,  $\lambda_1 \rightarrow 2 \sqrt{2}$ and $\lambda_2 \rightarrow 0 \,$.
\end{itemize}

\label{thm:stable}
\end{theorem}

\begin{figure}[tb]
\begin{center}
\includegraphics[width=240bp]{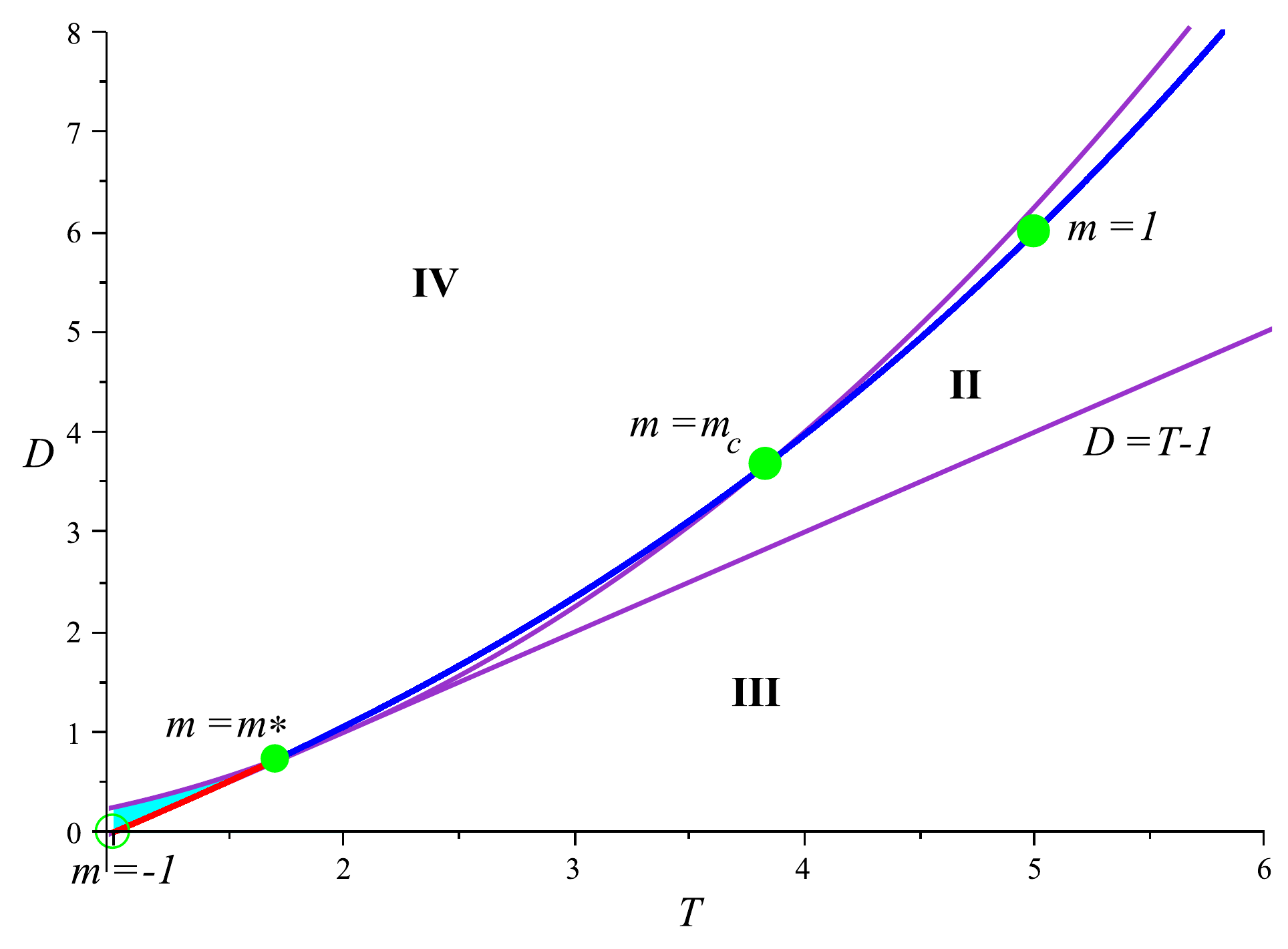}
\hspace{0.1in}
\includegraphics[width=240bp]{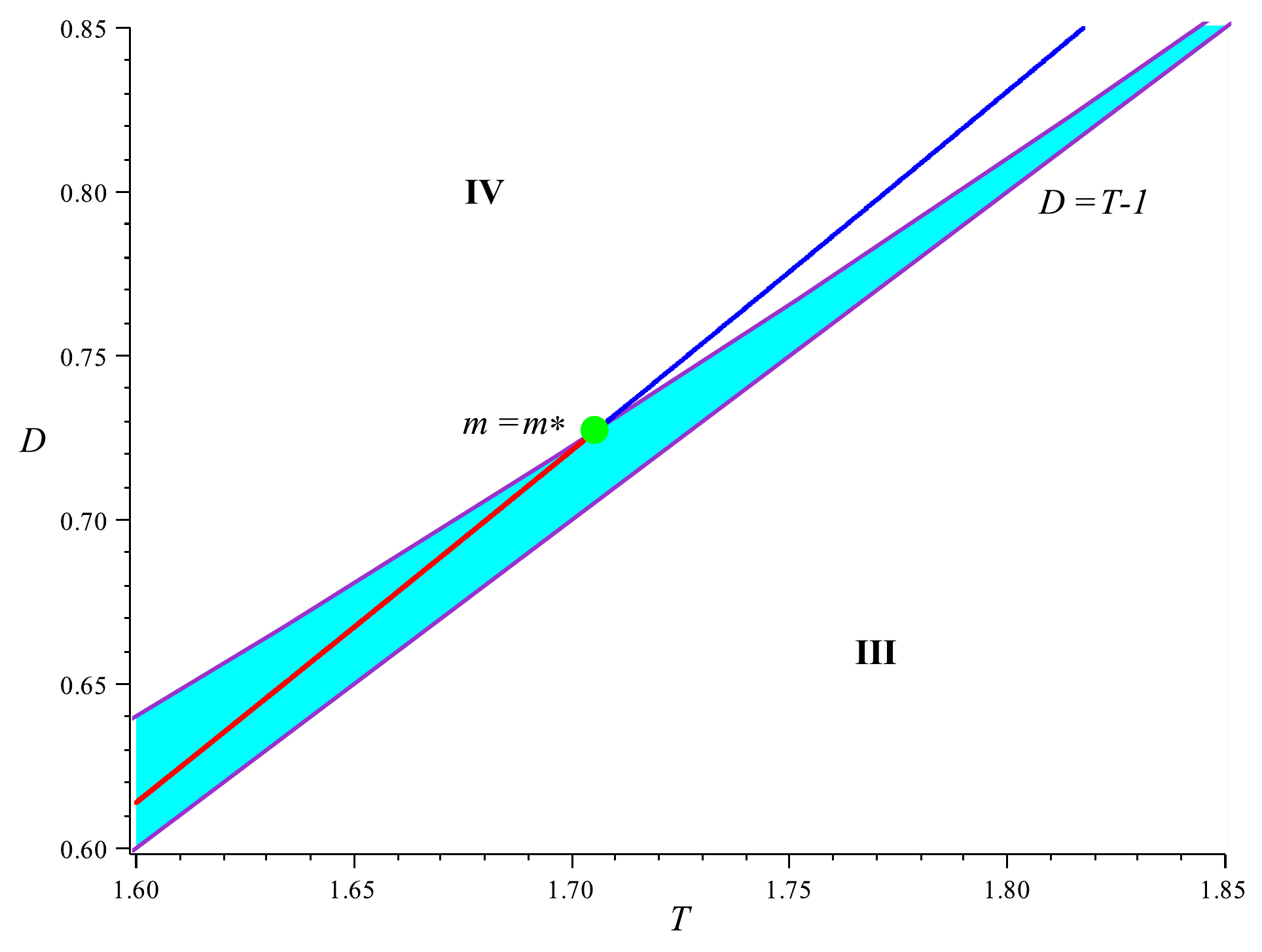}
\end{center}
\vspace{-0.2in}
\caption{The trace-determinant curve as $m$ varies for the Group~I solutions (red is stable, blue is unstable).  For $m > m_c$, the curve lies in region~II and the solutions have two real pairs of eigenvalues.
At $m = m_c \approx -0.0175$, the eigenvalues bifurcate into a complex quartuplet as the curve enters region~IV.  As $m$ decreases past $m^{\ast} \approx -0.8564$, the eigenvalues bifurcate again, forming two pairs of
pure imaginary eigenvalues.  The solutions are linearly stable for $-1 < m < m^{\ast}$.  A zoom of the figure near the key bifurcation value $m^{\ast}$ is shown to the right.}
\label{fig:TD-sols}
\end{figure}

\pf  We begin by focusing on the solutions with the orderings in Group~I.  By Theorem~\ref{thm:TDinv} we can restrict our attention to one particular solution from this group.
Since the solution varies continuously in $m$ ($m > -1$), so do the values of $T = T(m)$ and $D = D(m)$ that govern stability.  
Figure~\ref{fig:TD-sols} shows a plot of $D$ versus $T$ as $m$ varies, including two bifurcations at
$m = m_c$ and $m = m^{\ast}$.   Our intent is to justify this picture rigorously.

Using the asymptotic expansion for the solution with ordering $\ord{1}{3}{2}{4}$, we find from equations (\ref{eq:TrAs}) and~(\ref{eq:DtAs}) that
$T$ and $D$ are approaching $1$ and $0$, respectively, as $m \rightarrow -1^{+}$.   By formula~(\ref{eq:nontrivEvals}),
the nontrivial normalized eigenvalues are limiting on $0, 0, \pm i$.
We also find that 
$$
D - (T - 1) \; = \;  18 \epsilon^6  - 162 \epsilon^8 +  \mathcal{O}(\epsilon^{10}) ,
$$
so that $D > T - 1$ for $m$ sufficiently close to $-1$.  
Since the other three conditions in Theorem~\ref{thm:TDinv} are also satisfied, 
the solution is linearly stable.

This proves that the trace-determinant curve for the Group~I solutions lies in the stability region for $m$
sufficiently close to $-1$.  At $m = 1$, this curve reaches the point $(T=5, D=6)$ in the unstable region~II
(two real pairs of eigenvalues).  To ascertain how stability is lost, we search for bifurcations, that is, we
look for intersections between the trace-determinant curve and the boundaries of the stability region.

Adding the polynomial obtained from the numerator
of $T^2/4 - D$ to our defining system of equations $\{f_1, f_2, \ldots, f_6\}$ yields a system of polynomials whose solutions
contain those relative equilibria with repeated eigenvalues.  Fortunately, it is possible to compute a lex Gr\"{o}bner basis for this augmented system and
eliminate all variables except for the parameter~$m$.  
The first polynomial in this basis is $(m+1)^2 (m+3)^2 \Psi(m)$.  Therefore, repeated eigenvalues may only occur
if $m = m_c$ or $m = m_{\ast}$.

Applying back-substitution into the Gr\"{o}bner basis, we find that there are six solutions corresponding to the orderings in Group~I
at both $m = m_c$ and $m = m^{\ast}$.  At $m = m_c$, the solutions have a $T$-value larger than 2 ($T \approx 3.8344$) and are
therefore unstable.   On the other hand, the value of $T$ for the solutions at $m = m^{\ast}$ is $T^{\ast} \approx 1.7054$.  
Applying formula~(\ref{eq:nontrivEvals}), the
Group~I relative equilibria are spectrally stable with repeated nontrivial normalized eigenvalues $\pm i \sqrt{1 - (T^{\ast}/2)^2}$.
The relative equilibria at $m = m^{\ast}$ are not linearly stable because the matrix C is not a scalar multiple of the identity matrix.
This was confirmed by appending the numerator of $C_{21}$ to the augmented system described above and computing a Gr\"{o}bner basis.
Since the polynomial~1 was obtained, the value of $C_{21}$ at either bifurcation is nonzero.

Computing the trace and determinant for the Group~I solution at $m = -1/2$ ($\rho = -\sqrt{3/7}$), we find that $D > T^2/4$ and thus
the nontrivial eigenvalues $\mu_1$ and $\mu_2$ are complex.   It follows that the trace-determinant curve for the Group~I solution
lies in region~IV for $m^{\ast} < m < m_c$.   By part~{\bf (ii)} of Lemma~\ref{lemma:invsubspace}, the normalized
eigenvalues form a complex quartuplet $\pm \alpha \pm i \beta$ for these $m$-values.

Next, we add the numerator of $D - (T - 1)$ to the system $\{f_1, f_2, \ldots, f_6\}$ and compute a lex Gr\"{o}bner basis for this augmented system.
The first term in the Gr\"{o}bner basis is simply $(m+1)^2 (m+3)^2$.  The $m+1$ term is expected because $T=1$ and $D=0$ are the limiting
values as $m \rightarrow -1^{+}$.  The fact that there are no other roots for $m > -1$ shows that $D > T - 1$ for all solutions (using continuity) because 
the values at $m = 1$ satisfy this inequality (see Example~\ref{ExampleM1Stab}).  
This is true for solutions from either Group I or~II.  
A similar Gr\"{o}bner basis calculation shows that $T > 1$ for all solutions.  It follows that the trace-determinant curve
lies in the first quadrant and never in region~III.  Moreover, the only possible bifurcations occur at $m = m_c$ and $m = m^{\ast}$,
where the curve crosses the repeated root parabola $D = T^2/4$.  Thus, we have shown that the trace-determinant curve for the Group~I solutions lies in the stability region
for $-1 < m < m^{\ast}$, in region~IV for $m^{\ast} < m < m_c$, and in region~II for $m > m_c$.

To determine the fate of the Group~I normalized eigenvalues as $m \rightarrow \infty$, we compute an asymptotic expansion for the solution with ordering
$\ord{1}{2}{3}{4}$.  Setting $m = 1/\epsilon^2$ and treating $\epsilon$ as a small parameter, we find
\begin{eqnarray*}
x_3 & = &  3  - \frac{\sqrt{3}}{3} \epsilon + \frac{1}{12} \epsilon^2 + \frac{1025}{1296} \sqrt{3} \, \epsilon^3 - \frac{2059}{5184} \epsilon^4 +   \mathcal{O}(\epsilon^{5}), \\[0.07in]
x_4 & = &  \frac{4}{3} \sqrt{3} \, \epsilon^{-1} + \frac{2}{3}  + \frac{35}{27} \sqrt{3} \, \epsilon + \frac{307}{648} \epsilon^2 +   \mathcal{O}(\epsilon^{3}), \\[0.07in]
T & = &  8 - \frac{45}{4} \epsilon^2 +   \mathcal{O}(\epsilon^{4}),  \\[0.07in]
D & = &  15 - \frac{171}{4} \epsilon^2 +   \mathcal{O}(\epsilon^{4}) \, .
\end{eqnarray*}
Thus, $\mu_1 = (T + \sqrt{T^2 - 4D})/2$ approaches $5$ and $\mu_2 = (T - \sqrt{T^2 - 4D})/2$ approaches $3$ as $m \rightarrow \infty$.
The limiting behavior of the nontrivial normalized eigenvalues now follows from formula~(\ref{eq:nontrivEvals}).
We note that the solution with this particular ordering limits on a configuration with the three equal-strength vortices equally spaced
($r_{12} = r_{23} = 2$) and the fourth vortex infinitely far away.  This completes the proof of items {\bf (i)} and~{\bf (ii)} of the theorem.

The Gr\"{o}bner basis calculations above show that the Group~II solutions never bifurcate.  Since we also have $T = 5$ and $D = 6$ at $m=1$ for this
group of orderings, it follows that the trace-determinant curve for the Group~II solutions is always contained in region~II.  Consequently, the
nontrivial eigenvalues form two real pairs $\pm \lambda_1, \pm \lambda_2$ for all $m > -1/2$.  Using the asymptotic expansions (\ref{eq:x3Gr2As}) and~(\ref{eq:x4Gr2As})  
for the Group~II solution $\ord{1}{2}{4}{3}$, we find that 
$$
T \; = \;  \frac{21}{10} \epsilon^{-2} + \frac{578}{175} +  \mathcal{O}(\epsilon^{2})  \quad \mbox{and } \quad
D \; = \;  \frac{189}{50} \epsilon^{-2} + \frac{3459}{875} +  \mathcal{O}(\epsilon^{2})
$$
are expansions for the trace and determinant of the Group~II solutions for $m$ close to $-1/2$.
It follows that $\mu_1 \approx (21/10) \epsilon^{-2}$ approaches $\infty$, while $\mu_2$
approaches $9/5$ as $m \rightarrow -1/2^{+}$.  The limiting behavior
of the nontrivial normalized eigenvalues now follows from formula~(\ref{eq:nontrivEvals}).

To determine the fate of the Group~II normalized eigenvalues as $m \rightarrow \infty$, we compute an asymptotic expansion for the solution with ordering
$\ord{1}{4}{3}{2}$.  Setting $m = 1/\epsilon^2$ and treating $\epsilon$ as a small parameter, we find
\begin{eqnarray*}
x_3 & = &  1  -  \epsilon + \frac{1}{4} \epsilon^2 + \frac{1}{16} \epsilon^3 - \frac{3}{64} \epsilon^4 +   \mathcal{O}(\epsilon^{5}), \\[0.07in]
x_4 & = &  -\frac{1}{4} \epsilon - \frac{3}{8} \epsilon^2 + \frac{1}{8} \epsilon^3 +   \mathcal{O}(\epsilon^{4}), \\[0.07in]
T & = &  4 + \frac{3}{4} \epsilon^2 +   \mathcal{O}(\epsilon^{4}),  \\[0.07in]
D & = &  3 + \frac{21}{4} \epsilon^2 +   \mathcal{O}(\epsilon^{4}) \, .
\end{eqnarray*}
It follows that $\mu_1 \rightarrow 3$ and $\mu_2 \rightarrow 1^{+}$ as $m \rightarrow \infty$.
The limiting behavior of the nontrivial normalized eigenvalues now follows from formula~(\ref{eq:nontrivEvals}).
We note that the solution with this particular ordering limits on a configuration containing a collision between vortices
2 and~3, with the fourth vortex located in the middle of vortices 1 and~2 ($r_{14} = r_{24} = 1$).
This completes the proof of item~{\bf (iii)}.
\enpf

\begin{remark}

\begin{enumerate}

\item  Using Gr\"{o}bner bases, it is possible to express the $x_3$-coordinate at $m = m^{\ast}$ as the root 
of an even 36th-degree polynomial in one variable with integer coefficients.
The same is true for the $x_4$-coordinate (same degree, different polynomial).  

\item  The fact that both the Group I and~II collinear relative equilibria are unstable for $m > 0$ agrees with Corollary~3.5 in~\cite{g:stability}.
In general, any collinear relative equilibrium of $n$ vortices, where all circulations have the same sign, is always unstable and has
$n-2$ nontrivial real eigenvalue pairs $\pm \lambda_j$.  This follows by generalizing a clever argument of Conley's from the collinear $n$-body setting
(see Pacella~\cite{pacella} or Moeckel~\cite{rick-book} for details in the $n$-body case).

\item  Recall that relative equilibria are critical points of the Hamiltonian $H$ restricted to a level surface of the angular
impulse~$I$.   Numerical calculations in Matlab indicate that all of our solutions (both stable and unstable) are saddles
(the Morse index is always 2, except for $m = 0$).   Thus, in contrast to the case of same-signed circulations, with mixed signs it is possible for a saddle to be linearly stable.
A similar observation, using a modified potential function, was also made in~\cite{bhl}.

\end{enumerate}

\end{remark}

\section{Conclusion}

We have used ideas from modern and computational algebraic geometry to rigorously study the
collinear relative equilibria in the four-vortex problem where three circulation strengths are assumed identical.
Exploiting the $S_3$ invariance in the problem, we simplified the defining equations and obtained
a specific count on the number and type of solutions in terms of the fourth vorticity $\Gamma_4 = m$.  
The linear stability of solutions in the full plane was investigated and stable solutions were discovered 
for $m$ negative.   Reductions were made to simplify the stability calculations and useful formulas were
derived that apply to any four-vortex collinear relative equilibrium.  Asymptotic expansions
were computed to rigorously justify the behavior of solutions near collision.  Gr\"{o}bner bases were
used to locate key bifurcation values.
It is hoped that the reductions employed here involving symmetry and invariant group
theory will prove useful in similar problems.

\vs

\noindent  {\bf Acknowledgments:} 
The authors would like to thank the National Science Foundation (grant DMS-1211675) and the
Holy Cross Summer Research Program for their support.

\bibliographystyle{amsplain}

\begin{thebibliography}{99}

\bibitem{albouy}  Albouy, A., Fu, Y., Sun, S., Symmetry of planar four-body convex central configurations,
{\em Proc. R. Soc. Lond. Ser. A Math. Phys. Eng. Sci.}  {\bf 464} (2008), no. 2093, 1355--1365. 


\bibitem{aref-equil} Aref, H., On the equilibrium and stability of a row of point vortices,
{\em J. Fluid Mech.} {\bf 290} (1995), 167--181.

\bibitem{aref-int} Aref, H., Integrable, chaotic, and turbulent vortex motion in two-dimensional flows,
{\em Ann. Rev. Fluid. Mech.} {\bf 15} (1983), 345--389.




\bibitem{aref-stab3} Aref, H., Stability of relative equilibria of three vortices,
{\em Phys. Fludis} {\bf 21} (2009), 094101.


\bibitem{aref-newton} Aref, H., Newton, P. K., Stremler, M. A., Tokieda, T., Vainchtein, D. L., Vortex crystals,
{\em Adv. Appl. Mech.} {\bf 39} (2003), 1--79.


\bibitem{bhw} Barry, A. M., Hall, G. R., Wayne, C. E., 
Relative equilibria of the $(1+n)$-vortex problem,
{\em J. Nonlinear Sci.} {\bf 22} (2012), 63--83.


\bibitem{bhl} Barry, A. M., Hoyer-Leitzel, A., Existence, stability, and symmetry of relative equilibria with a dominant vortex,
{\em SIAM J. Appl. Dyn. Syst.} {\bf 15}, no. 4 (2016), 1783--1805.



\bibitem{cs} Cabral, H. E., Schmidt, D. S., Stability of relative equilibria in
the problem of $N+1$ vortices, {\em SIAM J. Math Anal.}
{\bf 31}, no. 2 (1999), 231--250.




\bibitem{CLO}  Cox, D. A., Little, J. B., O'Shea, D., {\em Ideals, Varieties, and Algorithms:
An Introduction to Computational Algebraic Geometry and Commutative Algebra}, 3rd ed.,
Springer, Berlin (2007). 


\bibitem{davis} Davis, C., Wang, W., Chen, S. S.,  Chen, Y., Corbosiero, K., DeMaria, M., 
Dudhia, J., Holland, G., Klemp, J., Michalakes, J., Reeves, H., Rotunno, R., Snyder, C.,
Xiao, Q.,  
Prediction of Landfalling Hurricanes with the Advanced Hurricane WRF Model,
{\em Monthly Weather Review} {\bf 136} (2007), 1990--2005.





\bibitem{HRS} Hampton, M., Roberts, G. E., Santoprete, M.,
Relative equilibria in the four-vortex problem with two pairs of equal vorticities,
{\em J. Nonlinear Sci.} {\bf 24} (2014), 39-92.


\bibitem{have} Havelock, T. H., The stability of motion of rectilinear vortices in
ring formation, {\em Philosophical Magazine} {\bf 11}, no. 7 (1931), 617--633.








\bibitem{kirchhoff}  Kirchhoff G., Vorlesungen \"{u}ber Mathematische Physik, {\bf I}, 
Teubner, Leipzig, 1876.


\bibitem{KossSchub} Kossin, J. P., Schubert, W. H., Mesovortices, polygonal flow patterns, 
and rapid pressure falls in hurricane-like vortices, {\em J. Atmos. Sci.}
{\bf 58} (2001), 2196--2209.



\bibitem{maple} Maple, version 15.00, (2011), Maplesoft, Waterloo Maple Inc.

\bibitem{matlab} MATLAB, version 7.10.0.499 (R2010a), (2010), The MathWorks, Inc.


\bibitem{meyer} Meyer, K. R., Hall, G. R., Offin, D., {\em Introduction to Hamiltonian Dynamical Systems
and the $N$-Body Problem}, 2nd ed., Applied Mathematical Sciences, 90, Springer, New York (2009).




\bibitem{rick-book}  Moeckel, R., Central configurations, in
{\em Central Configurations, Periodic Orbits, and Hamiltonian Systems},
Llibre, J., Moeckel, R., Sim\'{o}, C, Birkh\"{a}user (2015), 105--167.








\bibitem{moulton}  Moulton, F. R.,  The straight line solutions of the problem of $n$~bodies, 
{\em Ann. of Math. (2)} {\bf 12}, no. 1 (1910), 1--17. 


\bibitem{newton} Newton, P. K., {\em The $N$-Vortex Problem: Analytic Techniques}, 
Springer, New York (2001).





\bibitem{oneil}  O'Neil, K. A.,  Stationary configurations of point vortices, 
{\em Trans. Amer. Math. Soc.}  {\bf 302}, no. 2 (1987),  383--425. 


\bibitem{pacella}  Pacella, F., Central configurations of the $N$-body 
problem via equivariant Morse theory,  {\em Arch. Ration. Mech. Anal.}
{\bf 97} (1987), 59--74.


\bibitem{palmore}  Palmore, J.,  Relative equilibria of vortices in two dimensions,
{\em Proc. Natl. Acad. Sci. USA} {\bf 79} (Jan. 1982), 716--718. 


\bibitem{perez}  P\'{e}rez-Chavela, E., Santoprete, M., Tamayo, C., 
Symmetric relative equilibria in the four-vortex problem with three equal vorticities,
{\em Dyn. Contin. Discrete Impuls. Syst. Ser. A Math. Anal.}  {\bf 22}, no. 3 (2015), 189--209.



\bibitem{g:stability}  Roberts, G. E., Stability of relative equilibria in the planar $n$-vortex problem,
{\em SIAM J. Appl. Dyn. Syst.} {\bf 12}, no. 2 (2013), 1114--1134.


\bibitem{schm} Schmidt, D., The stability of the Thomson heptagon,
{\em Regul. Chaotic Dyn.} {\bf 9}, no. 4 (2004), 519--528.


\bibitem{spring} Spring, D., On the second derivative test for constrained local extrema,
{\em Amer. Math. Monthly} {\bf 92} (1985), no. 9, 631--643.
 
 
\bibitem{sturm}  Sturmfels, B., {\em Solving Systems of Polynomial Equations},
Conference Board of the Mathematical Sciences Regional Conference Series in Mathematics,
no. 97, Amer. Math. Soc.  (2002).  


\bibitem{synge} Synge, J. L., On the motion of three vortices, {\em Can. J. Math}
{\bf 1} (1949), 257--270.


\bibitem{thomson}  Thomson, J. J., {\em A Treatise on the Motion of Vortex Rings: An essay to which the Adams 
prize was adjudged in 1882}, University of Cambridge, Macmillan, London (1883).


\bibitem{xia} Xia, Z., {\em Central configurations with many small masses}, J. Differential Equations {\bf 91} (1991), 168--179.



\end{thebibliography}

\end{document}